\pdfoutput=1
\documentclass[10pt]{article}
\usepackage{graphicx}
\usepackage{amssymb}
\usepackage{epstopdf}
\usepackage{amsmath}
\usepackage{amsfonts}
\newcommand{\braket}[2]{\langle #1,#2 \rangle}

\def\phi{{\varphi}}

\DeclareSymbolFont{AMSb}{U}{msb}{m}{n}
\DeclareMathSymbol{\N}{\mathbin}{AMSb}{"4E}
\DeclareMathSymbol{\Z}{\mathbin}{AMSb}{"5A}
\DeclareMathSymbol{\R}{\mathbin}{AMSb}{"52}
\DeclareMathSymbol{\Q}{\mathbin}{AMSb}{"51}
\DeclareMathSymbol{\I}{\mathbin}{AMSb}{"49}
\DeclareMathSymbol{\C}{\mathbin}{AMSb}{"43}

\def\be{\begin{equation}}
\def\ee{\end{equation}}
\def\ber{\begin{eqnarray}}
\def\eer{\end{eqnarray}}

\def\beq{\begin{equation}}
\def\eeq{\end{equation}}

\begin{document}

\addtolength{\textheight}{0 cm} \addtolength{\hoffset}{0 cm}
\addtolength{\textwidth}{0 cm} \addtolength{\voffset}{0 cm}

\newcommand{\ZZ}{\mathbb{Z}}
\newcommand{\Rm}{\mathbb{R}}
\newcommand{\RR}{\mathbb{R}}
\newcommand{\NN}{\mathbb{N}}
\newcommand{\sU}{\mathcal{U}}
\newcommand{\sF}{\mathcal{F}}
\newcommand{\sM}{\mathcal{M}}
\newcommand{\sS}{\mathcal{S}}
\newcommand{\mL}{\mathcal{L}}
\newcommand{\ac}{\hbox{\small ac}}
\newcommand{\mC}{\ensuremath{\mathcal{C}}}
\newcommand{\mU}{\ensuremath{\mathcal{U}}}
\newcommand{\mT}{\ensuremath{\mathcal{T}}}
\newcommand{\mS}{\ensuremath{\mathcal{S}}}
\newcommand{\mF}{\ensuremath{\mathcal{F}}}
\newcommand{\Nm}{\ensuremath{\mathbb{N}}}
\newcommand{\Zm}{\ensuremath{\mathbb{Z}}}
\newcommand{\Hm}{\ensuremath{\mathbb{H}}}
\newcommand{\mM}{\ensuremath{\mathcal{M}}}
\newcommand{\mK}{\ensuremath{\mathcal{K}}}
\newcommand{\mD}{\ensuremath{\mathcal{D}}}
\newcommand{\mA}{\ensuremath{\mathcal{A}}}
\newcommand{\mO}{\ensuremath{\mathcal{O}}}
\newcommand{\mI}{\ensuremath{\mathcal{I}}}
\newcommand{\mB}{\ensuremath{\mathcal{B}}}
\newcommand{\Tm}{\ensuremath{\mathbb{T}}}
\newcommand{\mE}{\ensuremath{\mathcal{E}}}
\newcommand{\vs}{\vspace{.5cm}}

\def\proof {\noindent{\sc{Proof. }}}
\def\qed {\mbox{}\hfill {\small \fbox{}} \\}  
\def\lto{\longrightarrow}
\def\lmto{\longmapsto}
\def\eq{\Longleftrightarrow}
\def\leq{\leqslant}
\def\geq{\geqslant}

\newtheorem{lem}{Lemma}
\newtheorem{thm}{Theorem}
\newtheorem{conj}[lem]{Conjecture}
\newtheorem{ques}[lem]{Question}
\newtheorem{cor}[lem]{Corollary}
\newtheorem{prop}[lem]{Proposition}
\newtheorem{defn}[lem]{Definition}
\newtheorem{note}[lem]{Note}
\newtheorem{rmk}{Remark}
\def\proof {\noindent{\sc{Proof. }}}
\def\qed {\mbox{}\hfill {\small \fbox{}} \\}  

\newcommand{\grad}{\operatorname{grad}}
\newcommand{\Leg}{\mathcal{L}}

\title{Optimal Ballistic Transport and Hopf-Lax Formulae on Wasserstein Space}
\author{Nassif  Ghoussoub\thanks{Partially supported by a grant from the Natural Sciences and Engineering Research Council of Canada. This work was initiated during a visit of the author to the Schr\"odinger Institute in Vienna during the month of June 2016.}\\ \\
{\it\small Department of Mathematics,  University of British Columbia}\\
{\it\small Vancouver BC Canada V6T 1Z2}\\
{\small nassif@math.ubc.ca}\vspace{1mm}
}
\maketitle

\begin{abstract}  We investigate the optimal mass transport problem associated to the following ``ballistic"  cost functional on phase space $M\times M^*$,  
$$
b_T(v, x):=\inf\{\langle v, \gamma (0)\rangle +\int_0^TL(\gamma (t), {\dot \gamma}(t))\, dt; \gamma \in C^1([0, T), M);   \gamma(T)=x\}, 
$$
where $M=\R^d$, $T>0$, and $L:M\times M \to \R$ is a Lagrangian that is jointly convex in both variables. Under suitable conditions on the initial and final probability measures, we use convex duality \`a la Bolza and Monge-Kantorovich theory to lift classical Hopf-Lax formulae from state space to Wasserstein space.  This allows us to relate  
optimal transport maps for the ballistic cost 
to those associated with the fixed-end cost defined on $M\times M$ by
$$
c_T(x,y):=\inf\{\int_0^TL(\gamma(t), {\dot \gamma}(t))\, dt; \gamma\in C^1([0, T), M);  \gamma(0)=x, \gamma(T)=y\}.
$$
We also point to links with the theory of mean field games.  

\end{abstract}

\section{Introduction and main results} Given a cost functional $c(y, x)$ on some product measure space $X_0\times X_1$, and  two probability measures $\mu$ on $X_0$ and $\nu$ on  $X_1$, we consider the problem of optimizing the total cost of transport plans and its corresponding dual principle as formulated by Kantorovich
\begin{equation*}
\inf\big\{\int_{X_0\times X_1} c(y, x)) \, d\pi; \pi\in \mK(\mu,\nu)\big\}=\sup\big\{\int_{X_1}\phi_1(x)\, d\nu(x)-\int_{X_0}\phi_0(y)\, d\mu(y);\,  \phi_1, \phi_0 \in \mK(c)\big\},
\end{equation*}
where $\mK(\mu,\nu)$ is the set of probability measures $\pi$ on $X_0\times X_1$ whose marginal on $X_0$ (resp. on $X_1$) is $\mu$ (resp., $\nu$) {\it (the transport plans)}, and where $\mK(c)$ is the set of functions $\phi_1\in L^1(X_1, \nu)$ and $\phi_0\in L^1(X_0, \mu)$ such that 
\[
\phi_1(x)-\phi_0(y) \leq c(y,x) \quad \hbox{for  all $(y,x)\in X_0\times X_1$}.
\]
The pairs of functions in $\mK(c)$ can be assumed to satisfy 
\[
\phi_1(x)=\inf_{y\in X_0} c(y, x)+\phi_0(y) \quad {\rm and} \quad\phi_0(y)=\sup_{x\in X_1} \phi_1(x)-c(y, x).
\]
They will be called {\it admissible Kantorovich potentials}, and for reasons that will become clear later, we shall say that $\phi_0$ (resp., $\phi_1$) is an initial (resp., final) Kantorovich potential.

The original Monge problem dealt with the cost $c(y, x)=|x-y|$ (\cite{M}, \cite{S}, \cite{E-G}, \cite{V1}, \cite{V2}) 
and was constrained to those probabilities in ${\mathcal K}(\mu, \nu)$ that are supported by graphs of measurable maps from $X$ to $Y$ pushing $\mu$ onto $\nu$.  Brenier \cite{B1} considered the important quadratic case  $c(x,y)=|x-y|^2$. This was followed by a large number of results addressing costs of the form $f(x-y)$, where $f$ is either a convex or a concave function \cite{G-M}. With a purpose of connecting mass transport with Mather theory, Bernard and Buffoni \cite{B-B} considered dynamic cost functions on a given compact manifold $M$, that deal with fixed end-points problems of the following type:
\begin{equation}\label{BB}
c_T(y,x):=\inf\{\int_0^TL(t, \gamma(t), {\dot \gamma}(t))\, dt; \gamma\in C^1([0, T), M);  \gamma(0)=y, \gamma(T)=x\},
\end{equation}
where $[0, T]$ is a fixed time interval, and $L: TM \to \R\cup\{+\infty\}$ is a given Lagrangian that is convex in the second variable of the tangent bundle $TM$. Fathi and Figalli \cite{F-F} eventually dealt with the case where $M$ is a non-compact Finsler manifold. Note that standard cost functionals of the form  $f(|x-y|)$, where $f$ is convex, are particular cases of the dynamic formulation, since they correspond to Lagrangians of the form $L(t, x, p)=f(p)$. 

In this paper, we shall consider the {\it ``ballistic cost function,"} which is defined on phase space $M^*\times M$ by, 
 \begin{equation}\label{bal}
b_T(v, x):=\inf\{\langle v, \gamma (0)\rangle +\int_0^TL(t, \gamma (t), {\dot \gamma}(t))\, dt; \gamma \in C^1([0, T), M);   \gamma(T)=x\}, 
\end{equation}
where $M$ is  a Banach space and $M^*$ is its dual. The associated transport problems will be 
\begin{equation}\label{bal1}
{\overline B}_T(\mu_0, \nu_T):=\sup \{\int_{M^*\times M} b_T(v, x)\, d\pi;\, \pi\in \mK(\mu_0,\nu_T)\}, 
\end{equation}
and
\begin{equation}\label{bal2}
{\underline B}_T(\mu_0, \nu_T):=\inf \{\int_{M^*\times M} b_T(v, x)\, d\pi;\, \pi\in \mK(\mu_0,\nu_T)\}, 
\end{equation}
where $\mu_0$ (resp., $\nu_T$) is a given probability measure on $M^*$ (resp., $M$). Note that when $T=0$, we have $b_0(x, v)=\langle v, x\rangle$, which is exactly the case considered by Brenier \cite{B1}, that is 
\begin{equation}
{\overline W}(\mu_0, \nu_0):=\sup \{\int_{M^*\times M} \langle v, x\rangle\, d\pi;\, \pi\in \mK(\mu_0,\nu_0)\}, 
\end{equation}
and 
\begin{equation}
{\underline W}(\mu_0, \nu_0):=\inf \{\int_{M^*\times M} \langle v, x\rangle\, d\pi;\, \pi\in \mK(\mu_0,\nu_0)\}, 
\end{equation}
making (\ref{bal1}) a suitable dynamic version of the Wasserstein distance.

The assumptions on the Lagrangian that we use are as follows:
\begin{trivlist}

\item (A1) $L:M\times M \to \R\cup\{+\infty\}$ is convex, proper and lower semi-continuous in both variables.

\item (A2) The set $F(x):=\{p; L(x, p) <\infty\}$ is non-empty for all $x\in M$, and for some $\rho>0$, we have ${\rm dist} (0, F(x)) \leq \rho (1+|x|)$ for all $x\in M$.

\item (A3) For all $(x, p)\in M\times M$, we have $L(x,p) \geq \theta (\max\{0, |p|-\alpha |x|\})-\beta |x|$, where $\alpha, \beta$ are constants, and $\theta$ is a coercive, proper, non-decreasing function on $[0, \infty)$. 
\end{trivlist}
The associated Hamiltonian on $[0, T] \times M\times M^*$ is defined by
\[
H(t, x, q)=\sup_{p\in M}\{\langle p, q\rangle -L(t, x, p)\}.
\]
We shall assume throughout that $M=M^*=\R^d$, while preserving --for pedagogical reasons-- the notational distinction between the state space and its dual. These conditions on the Lagrangian make sure that the Hamiltonian $H$ is finite, concave in $x$ and convex in $q$, hence locally Lipschitz. Moreover, we have
\begin{equation}
\psi(x)-(\gamma |x| +\delta)|q| \leq H(x, q) \leq \phi(q)+(\alpha |q| +\beta)|x| \hbox{for all $x, q$ in $M\times M^*$,}
\end{equation}
where $\alpha, \beta, \gamma, \delta$ are constants, $\phi$ is finite and convex and $\psi$ is finite and concave (see \cite{R-W2}.\\
We note that under these conditions, the cost $(x, y)\to c(t, x,y)$ is convex proper and lower semi-continuous on $M\times M$. But the cost $b_T$ is nicer in many ways. For one, it is everywhere finite and locally Lipschitz continuous on $[0, \infty)\times M\times M^*$. Moreover, 
we shall consider suitable solutions of the Hamilton-Jacobi equation, 
\begin{eqnarray}\label{HJ.0} 
\left\{ \begin{array}{lll}
\partial_tV+H(t, x, \nabla_xV)&=&0 \,\, {\rm on}\,\,  [0, T]\times M,\\
\hfill V(0, x)&=&V_0(x), 
\end{array}  \right.
 \end{eqnarray}
which are formally given by the formula 
\begin{equation}\label{value.1}
V(t,x)=\inf\Big\{V_0(\gamma (0))+\int_0^tL(s,\gamma (s), {\dot \gamma}(s))\, ds; \gamma \in C^1([0, T), M);   \gamma(t)=x\Big\},
\end{equation}
as well as the following  dual Hamilton-Jacobi equation:
\begin{eqnarray}\label{dHJ} 
\left\{ \begin{array}{lll}
\partial_tW-H(\nabla_vW, v)&=&0 \,\, {\rm on}\,\, [0, T]\times M^*,\\
\hfill W(T, v)&=&W_T(v),  
\end{array}\right.
\end{eqnarray} 
whose variational solution is given by
\begin{equation}\label{value.2}
W(t,v)=\sup\Big\{W_T(\gamma (T))-\int_0^t{\tilde L}(s,  \gamma (s), {\dot \gamma}(s))\, ds; \gamma \in C^1([0, T), M^*);   \gamma(0)=v\Big\},
\end{equation}
where the Lagrangian ${\tilde L}$ is defined on $M^*\times M^*$ by 
$$
\tilde L(t, v,q):=L^*(t, q, v)=\sup\{\langle v, y\rangle +\langle p, q\rangle -L(t, y, p);\, (y, p)\in M\times M\}.
$$
While the above formula for the ``solutions" $V$ (resp., $W$) of (\ref{HJ.0}) (resp., (\ref{dHJ}))  is natural, it doesn't often generate solutions with suitable regularity properties, unless more conditions are imposed on the Lagrangian.  Bernard and Buffoni \cite{B-B} imposed conditions  that allowed for the consideration of so-called continuous  {\it viscosity solutions}. There are however instances where our convex setting is completely satisfactory under the above assumptions on $L$. Indeed, if $V_0$ (resp., $W_T$) is convex (resp., concave), then the solution $V(t, \cdot)$ (resp.,  $W(t, \cdot)$) is convex (resp., concave) for each $t$, and convex subdifferentiability then provides a good alternative to regularity. This is the case for $b(t, v, x)$, which remarkably satisfies both equations (in the sense of convex analysis), that is 
\begin{eqnarray*}\label{HJ}
\partial_tb+H(t, x, \partial_xb)&=&0 \quad {\rm on} \quad [0, T]\times M,\quad b(0, v, x)=\langle v, x\rangle.\\
\partial_tb-H(t, {\tilde \partial}_v b, v)&=&0 \quad {\rm on} \quad [0, T]\times M^*,\quad b(T, v, x)=\langle v, x\rangle.
\end{eqnarray*} 
Unfortunately, as we shall see below, the mass transport problems that we consider lead to an opposite situation, where the initial  function $V_0$ is concave, a property that is not propagated forward by $L$, and the terminal function $W_T$ is convex, which is not propagated backward by ${\tilde L}$. We shall therefore only consider {\it variational solutions} for (\ref{HJ.0}) (resp., \ref{dHJ}) meaning those that has the form (\ref{value.1}) and (\ref{value.2}) respectively. 

The following duality will be proved in Section 5. 

\begin{thm} \label{duality} Assume $M=\R^d$ and that $L$ satisfies hypothesis (A1), (A2) and (A3), and let $\mu_0$ (resp. $\nu_T)$ be a probability measure on $M^*$ (resp., $M$). 
The following duality holds: 
\begin{enumerate}
\item Assume $\mu_0$ is absolutely continuous with respect to Lebesgue measure, then  
 \begin{equation}\label{ballistic.dual1}
{\underline B}_T(\mu_0,\nu_T)=\sup\left\{\int_MV(T,x)\, d\nu_T(x)+ \int_{M^*} {\tilde V_0}(v)\, d\mu_0(v); \,  \hbox{$V_0$ concave on $M$ \& $V$ solution of (\ref{HJ.0})} \right\}, 
\end{equation}
where ${\tilde V}_0$ is the concave Legendre transform of  $V_0$, i.e., 
\[
{\tilde V_0}(v)=\inf\{\langle v, y\rangle - V_0(y); \, y\in M\}.
\]
\item Assume $\nu_T$ is absolutely continuous with respect to Lebesgue measure, then
\begin{equation}\label{ballistic.dual2}
{\overline B}_T(\mu_0,\nu_T)=\inf\left\{\int_MW_T^*(x)\, d\nu_T(x)+ \int_{M^*} W (0, v)\, d\mu_0(v); \,  \hbox{$W_T$ convex on $M^*$ \& $W$ solution of (\ref{dHJ})} \right\}, 
\end{equation}
where $W_T^*$ is the convex Legendre transform of  $W_T$, i.e., 
$$
W^*_T(x)=\sup\{\langle v, x\rangle - W_T (v);\, v\in M^*\}.
$$
\end{enumerate}
\end{thm}
In order to investigate the support of optimal transport plans,  we shall need a subgradient form of Hamiltonian dynamics, a solution of which over a time interval  $[0, T]$ is any pair of $C^1$ arcs $(x(t), v(t))$ such that for almost evert $t\in 0, T]$, 
\begin{eqnarray}
\left\{ \begin{array}{lll}
\, \, \, \dot x(t)&\in& \partial_v H(x(t), v(t))\\
-\dot v(t)&\in& {\tilde \partial}_x H(x(t), v(t)). 
\end{array}\right.
\end{eqnarray} 
 The associated Hamiltonian flow is the one-parameter family of, possibly set-valued, mappings $(\phi^H_t)_t$,  
 \begin{equation}
\phi^H_t(x_0, v_0)=\{(x,v); \hbox{$\exists$ a Hamiltonian trajectory starting at $(x_0, v_0)$ with $(x(t),v(t))=(x,v)$\}.}
 \end{equation}
We shall prove the following:

\begin{thm} \label{attain} Assume $L$ satisfies hypothesis (A1), (A2) and (A3), and let $\mu_0$ (resp. $\nu_T)$ be a probability measure on $M^*$ (resp., $M$). 
\begin{enumerate}
\item If $\mu_0$ is absolutely continuous with respect to Lebesgue measure, then there exists a probability measure $\pi_0$ on $M^*\times M$, and a concave function $k: M \to \R$ such that   
\begin{equation}
{\underline B}_T(\mu_0,\nu_T)=\int _{M^*} b_T \large(v, x) d\pi_0,
\end{equation}
and $\pi_0$ is supported on the possibly set-valued map $v\to \pi^*\phi^H_T(\nabla {\tilde k}(v), v)$, 
with $\pi^*:M\times M^*\to M$ being the canonical projection, and $\tilde k$ is the concave Legendre transform of $k$. 
\item Assume $\nu_T$ is absolutely continuous with respect to Lebesgue measure, then there exists a probability measure ${\tilde \pi}_0$ on $M^*\times M$, and a convex function $h: M^* \to \R$ such that 
\begin{equation}
{\overline B}_T(\mu_0,\nu_T)=\int _{M^*} b_T \large(v, x) d{\tilde \pi}_0,
\end{equation}
and ${\tilde \pi_0}$ is supported on the possibly set-valued map 
$x\to \pi^*\phi^{-H_*}_T(\nabla h (x), x)$, where $\phi^{-H_*}_t$ is the Hamiltonian flow  associated to the Lagrangian $L_*(v,q)=L^*(-q, v)$, namely $H_*(x, q)=-H(-x, q)$.
\end{enumerate}
\end{thm}
If we assume further that $L$ is a Tonelli Lagrangian on $M\times M^*$, then the above maps are single-valued and completely solve the Monge version of the mass transport problems, that is 
\begin{equation}
{\underline B}_T(\mu_0,\nu_T)=\int _{M^*} b_T \large(v, \pi^*\phi^H_T(\nabla {\tilde k}(v), v)\large) d\mu_0(v), 
\end{equation}
and 
\begin{equation}
{\overline B}_T(\mu_0,\nu_T)= \int _{M} b_T \large(\pi^*\phi^{-H}_T(\nabla h^* (x), x), x\large) d\nu_T(x). 
\end{equation}
The above two theorems will follow from the following interpolation result, and what is known about 
the optimal mass transport  
\begin{equation}\label{BBT}
C_T(\nu_0, \nu_T):=\inf \{\int_{M\times M} c_T(x,y)\, d\pi;\, \pi\in \mK(\nu_0,\nu_T)\}, 
\end{equation}
where $\nu_0$ and $\nu_T$ are two given probability measures on $M$. We shall also need another cost function on $M^*\times M^*$, 
\begin{equation}\label{tilde}
{\tilde c}_T(u, v):=\inf\{\int_0^T{\tilde L}(t, \gamma (t), {\dot \gamma}(t))\, dt; \gamma\in C^1([0, T), M^*);  \gamma (0)=u, \gamma (T)=v\},
\end{equation}
and its associated transport
\begin{equation}\label{BBT*}
{\tilde C}_T(\mu_0, \mu_T):=\inf \{\int_{M^*\times M^*} {\tilde c}_T(x,y)\, d\pi;\, \pi\in \mK(\mu_0,\mu_T)\}. 
\end{equation}
\begin{thm} \label{interpol} Assume that $L$ satisfies hypothesis (A1), (A2) and (A3), and let $\mu_0$ (resp. $\nu_T)$ be a probability measure on $M^*$ (resp., $M$).
\begin{enumerate}

\item If $\mu_0$ is absolutely continuous with respect to Lebesgue measure, then   
\begin{equation}\label{HL.one}
{\underline B}_T(\mu_0,\nu_T)=\inf\{{\underline W}(\mu_0, \nu)+ C_T(\nu, \nu_T);\, \nu\in {\mathcal P}(M)\}.
\end{equation}
The infimum is attained at some probability measure $\nu_0$ on $M$, and the initial Kantorovich potential for $C_T(\nu_0, \nu_T)$ is concave. 
\item If $\nu_T$ is absolutely continuous with respect to Lebesgue measure, then  
\begin{equation}\label{HL.two}
{\overline B}_T(\mu_0,\nu_T)=\sup\{ {\overline W}(\nu_T, \mu)-{\tilde C}_T(\mu_0, \mu);\, \mu\in {\mathcal P}(M^*)\}.  
\end{equation}
The supremum is attained at some probability measure $\mu_T$ on $M^*$, and the final Kantorovich potential for ${\tilde C}_T(\mu_0, \mu_T)$ is convex.

\end{enumerate}
\end{thm}
The above interpolation formulas can be seen as extensions of those by Hopf-Lax on state space to Wassertsein space. Indeed, for any convex and lower semi-continuous function $g$, the associated value function
 \begin{eqnarray*}
\left\{ \begin{array}{lll}
V_g(t,x)&=&\inf\Big\{g(\gamma (0))+\int_0^tL(s,\gamma (s), {\dot \gamma}(s))\, ds; \gamma \in C^1([0, T), M);   \gamma(t)=x\Big\},\\
V(0, x)&=&g(x),
\end{array}\right.
\end{eqnarray*} 
can be written as 
\begin{equation}\label{HL}
V_g(t,x)=\inf\{g(y)+c(t, y, x);\, y\in M\}\quad \hbox{as well as \quad $V_g(t,x)=\sup\{b(t, v, x)-g^*(v); \, v\in M^*\}.$}
\end{equation}
In the case where the Lagrangian $L(x, p)=L_0(p)$ is only a function of $p$, and if $H_0$ is the associated Hamiltonian, then 
 \begin{equation}
 c_t(y, x)=tL_0(\frac{1}{t}|x-y|)\hbox{\quad  and \quad $b_t(v, x)=\langle v, x\rangle-tH_0(v),$}
 \end{equation}
 and (\ref{HL}) is nothing but the Hopf-Lax formula used to generate solutions for corresponding Hamilton-Jacobi equations.
Moreover, when $g$ is the linear functional $g(x)=\langle v, x\rangle$, then $b(t, v, x)$ is itself a solution to two Hamilton-Jacobi equations, since 
\begin{equation}\label{basic}
b(t, v, x)=\inf\{\langle v, y\rangle+c(t, y, x);\, y\in M\}=\sup\{\langle w,x\rangle-{\tilde c}(t, v,w);\, w\in M^*\}.
\end{equation}
Note that formulas (\ref{HL.one}) and (\ref{HL.two}) can now be seen as extensions of (\ref{basic}) to the space of probability measures, where the Wasserstein distance fill the role of the scalar product. 

Surprisingly, the extension of the dual formula 
\begin{equation}\label{babydual}
c(t,y,x)=\sup\{ b(t, v, x)-\langle v, y\rangle; v\in M^*\}.
\end{equation}
is not always possible. Indeed, we shall see that the key to proving (\ref{HL.one}) and (\ref{HL.two}) is that 
the initial Kantorovich potential for ${\underline B}_T(\mu_0, \nu_T)$ can be taken to be convex, while for ${\overline B}_T(\mu_0, \nu_T)$, the final Kantorovich potential can be assumed to be concave. This contrasts the case of the fixed enpoints cost $c_T$, which even though it is jointly convex in both variables, one cannot deduce much in terms of the convexity or concavity of the Kantorovich potentials corresponding to $C_T(\nu_0, \nu_T)$. Indeed, if $L(x, v)=\frac{1}{2}|v|^2$, which corresponds to the cost $c(y,x)=\frac{1}{2}|x-y|^2$, the initial Kantorovich potential is then of the form $\phi_0(y)=g(y)-\frac{1}{2}|y|^2$, where $g$ is a convex function. 

Note however, that the optimal transport $C_T(\nu_0, \nu_T)$ corresponding to the initial measure $\nu_0$ obtained via the factorization of ballistic optimal transport problems, has a concave initial Kantorovich potential. This turned out to be a necessary and sufficient condition. In the case where $L(x, v)=\frac{1}{2}|v|^2$, it says that for such a measure the initial Kantorovich potential $\phi_0$ can still be concave, that is 
$0\leq D^2g \leq I$. 

\begin{thm}\label{endpoints} Assume $M=\R^d$ and that $L$ satisfies hypothesis (A1), (A2) and (A3). Assume $\nu_0$ and $\nu_T$ are probability measures on $M$ such that $\nu_0$ is absolutely continuous with respect to Lebesgue measure. Then, the following are equivalent:

\begin{enumerate}

\item The initial Kantorovich potential of $C_T(\nu_0, \nu_T)$ is concave.
\item The following formula holds 
\begin{equation}\label{three}
C_T(\nu_0, \nu_T)=\sup\{{\underline B}_T(\mu,\nu_T)-{\underline W}(\nu_0, \mu);\, \mu\in {\mathcal P}(M^*)\}.  
\end{equation}
\end{enumerate}
\end{thm}
Here is an application.

\begin{cor} \label{cor1}Consider the cost $c(y, x)=c_0(x-y)$, where $c_0$ is a convex function on $M$ and let $\nu_0, \nu_1$ be  probability measures on $M$ such that the initial Kantorovich potential associated to $C_T(\nu_0, \nu_T)$ is concave. Then, there exist concave functions $\phi: M\to \R$ and $\psi: M^*\to \R$ such that
\begin{equation}
C_1(\nu_0, \nu_1)-K=\int _{M} c_0 ( \nabla \psi\circ \nabla \phi(y)-y) d\nu_0(y)=\int_M \langle \nabla {\tilde\psi}(y)-\nabla \phi (y), y\rangle \, d\nu_0(y), 
\end{equation}
where $K$ is a constant and ${\tilde \psi}$ is the concave Legendre transform of $\psi$.
\end{cor}
The interpolation formula can be seen as a Hopf-Lax formula on Wasserstein space, since for a fixed $\mu_0$ on $M^*$ (resp., fixed  $\nu_T$ on $M$), then as a function of the terminal (resp., initial) measure, we have 
 \begin{equation}
{\underline {\mathcal B}}^{\mu_0}(t, \nu)=\inf\{{\underline {\mathcal U}}^{\mu_0}(\rho)+ C_t(\rho, \nu);\, \rho \in {\mathcal P}(M)\}\,\, \hbox{and \,\, ${\overline {\mathcal B}}^{\nu_T}(t, \mu)=\inf\{{\overline {\mathcal U}}^{\nu_T}(\rho)- {\tilde C}_t(\rho, \mu);\, \rho\in {\mathcal P}(M^*)\},$ }
\end{equation}
where 
$${\underline {\mathcal U}}^{\mu_0}(\rho)={\underline W}(\mu_0, \rho)\hbox{\quad and \quad ${\overline {\mathcal U}}^{\nu_T}(\rho)={\overline W}(\nu_T, \rho).$}
$$  
 The following Eulerian formulation illustrates best how ${\underline {\mathcal B}}^{\mu_0}(t, \nu)$ and ${\overline {\mathcal B}}^{\nu_T}(t, \mu)$
 can be represented as value functionals on Wasserstein space. Indeed,  lift the Lagrangian $L$ to the tangent bundle of Wasserstein space via the formula 
\[
\hbox{${\mathcal L}(\rho, w);=\int_M L(x, w(x)) \, d\rho(x)$\quad and \quad $\tilde {\mathcal L}(\rho, w);=\int_{M^*} {\tilde L}(x, w(x)) \, d\rho(x),$
}
\]
where $\rho$ is any probability density  on $M$ (resp., $M^*$) and $w$ is a vector field on $M$ (resp., $M^*$). 
\begin{thm} Assume $L$ satisfies hypothesis (A1), (A2) and (A3), and let $\mu_0$ (resp. $\nu_T)$ be a probability measure on $M^*$ (resp., $M$). 
\begin{enumerate}
\item If $\mu_0$ is absolutely continuous with respect to Lebesgue measure, then  
\begin{eqnarray}
{\underline {\mathcal B}}^{\mu_0}(T, \nu):={\underline B}_T(\mu_0,\nu)
&=& \inf\left\{{\underline {\mathcal U}}^{\mu_0}(\rho_0) +
\int_0^T {\mathcal L}(\rho_t, w_t) dt;\,  \partial_t \varrho+ \nabla \cdot (\varrho w)=0,\,  \varrho_T=\nu\right\},  
 \end{eqnarray}
 
\item  If $\nu_T$ is absolutely continuous with respect to Lebesgue measure, then   
\begin{eqnarray}
{\overline {\mathcal B}}^{\nu_T}(T, \mu):={\overline B}_T(\mu, \nu_T)
&=& \sup\left\{ {\overline {\mathcal U}}^{\nu_T}(\rho_T) -
\int_0^T {\tilde {\mathcal L}}(\rho_t, w_t) dt; \, \,  \partial_t \varrho+ \nabla \cdot (\varrho w)=0,\,  \varrho_0=\mu\right\}. 
 \end{eqnarray}
The set of pairs $(\varrho, w)$ considered above are such that $t \rightarrow \varrho_t \in \mathcal P(M)$, and $t \rightarrow w_t (x)\in {\rm Lip}(\R^n)$ are paths of Borel vector fields.
\end{enumerate}
\end{thm}
One can then ask whether these value functionals also satisfy a Hamilton-Jacobi equation on Wasserstein space such as
  \begin{equation}\label{e:gangbo:master}
\left\{
\begin{array}{ll}
& \partial_t B 
+ {\mathcal H}(t, \nu, \nabla_\nu B(t, \nu)) 
=0,\\ [5pt] 
& B(0, \nu) ={\underline W}(\mu_0, \nu). 
\end{array}
\right.
\end{equation}  
Here the Hamiltonian  is defined as 
\[
{\mathcal H}(\nu, \zeta) =\sup\{\int \langle \zeta, \xi\rangle d\nu -{\mathcal L}(\nu, \xi); \xi \in T_\nu^*({\mathcal P}(M))\}.
\]
We note that Ambrosio-Feng \cite{A-F} have shown recently that --at least in the case where the Hamiltonian is the square-- value functionals on Wasserstein space yield a unique {\it metric viscosity solution} for (\ref{e:gangbo:master}). As importantly, Gangbo-Sweich \cite{G-S} have shown recently that under certain conditions, value functionals yield solutions to the so-called {\it Master equations} of mean field games. 

\begin{thm} (Gangbo-Swiech) Assume $\mathcal U_0: \mathcal P(M) \rightarrow \mathbb R$, and $U_0:  M \times \mathcal P(M) \rightarrow \mathbb R$ are functionals such that 
$ \nabla_x U_0(x, \mu) \equiv \nabla_\mu \mathcal U_0(\mu)(x)$ for all $x \in M$, $\mu \in \mathcal P(M),$ and consider the value functional,
\begin{equation*} 
\mathcal U(t, \nu)=\inf \left\{ \mathcal U_0(\varrho_0) + \int_0^t \mathcal L(\varrho, w) dt;\,  \partial_t \varrho+ \nabla \cdot (\varrho w)=0,\,  \varrho_T=\nu\right\}.
 \end{equation*} 
Then, there exists $U: [0,T] \times M \times \mathcal P(M) \rightarrow \mathbb R$ such that 
\begin{equation*}\label{1b}
 \nabla_x U_t(x, \nu) \equiv \nabla_\nu \mathcal U_t(\nu)(x) \quad \hbox{ for all $x \in M$, $\nu \in \mathcal P(M),$}
\end{equation*}  
 and  $U$ satisfies the Master equation below (\ref{e:gangbo:master.1}).
\end{thm} 
Applied to the value functional ${\underline {\mathcal B}}^{\mu_0}(t, \nu):={\underline B}_t(\mu_0,\nu)$, this 
should then yields the existence for any probabilities $\mu_0, \nu_T$, a function $\beta: [0,T] \times M \times \mathcal P(M) \rightarrow \mathbb R$ such that 
\begin{equation*}\label{1b}
 \nabla_x \beta (t, x, \nu) \equiv \nabla_\nu {\underline {\mathcal B}}^{\mu_0}(t, \nu)(x) \quad \hbox{ for all $x \in M$, $\nu \in \mathcal P(M),$}
\end{equation*}  
and  $\rho \in AC^2((0, T)\times {\mathcal P}(M))$ such that 
\begin{equation}\label{e:gangbo:master.1}
\left\{
\begin{array}{ll}
& \partial_t \beta  + \int\langle \nabla_\nu \beta (t, x, \nu)\cdot \nabla  H(x, \nabla_x\beta)\rangle\, d\nu
+ H(x, \nabla_x \beta(t, x,\nu)) =0,\\ [5pt] 
&\partial_t\rho +\nabla (\rho \nabla H(x, \nabla_x \beta))=0, \\ [5pt]
& \beta(0, \cdot, \cdot) = \beta_0, \quad \rho (T, \cdot)=\nu_T,
\end{array}
\right.
\end{equation}  
where $\beta_0(x, \rho)=\phi_\rho(x)$, where $\phi_\rho$ is the convex function such that $\nabla \phi_\rho$ pushes $\mu_0$ into $\rho$.  

Needless to say, one would like to consider value functionals on Wasserstein space that are more general than those starting with the Wasserstein distance. One can still obtain such functionals via mass transport by considering more general ballistic costs of the form 
\begin{equation}
b_g(T, v, x):=\inf\left\{g(v, \gamma (0)) +\int_0^TL(\gamma (t), {\dot \gamma}(t))\, dt; \gamma \in C^1([0, T), M)\right\},    
\end{equation}
where $g: M^*\times M \to \R$ is a suitable function. This, as well as the stochastic counterpart of the ballistic problem and the corresponding ``diffusive" master equation, will be considered  in a subsequent paper.

I am very grateful to Yann Brenier and Wilfrid Gangbo for several fruitful discussions, and in particular for pointing me towards the connection between value functionals on Wasserstein space and mean field games. I am also indebted to the Schr\"odinger Institute in Vienna, where most of this work was done during my visit there in June 2016.

\section{The Bolza duality and its consequences}

We first review properties of the various cost functions that will be needed in the sequel. They are standard and had been studied in detail in various articles by T. Rockafellar \cite{R1} and his co-authors \cite{R-W1, R-W2}. 

\begin{prop} Under assumptions (A) on the Lagrangian $L$, the costs $c$ and $b$ have the following properties:
\begin{enumerate}

\item For each $t\geq 0$, $(x, y)\to c(t, x,y)$ is convex proper and lower semi-continuous on $M\times M$.

\item For each $t\geq 0$, $v\to b(t, v, x)$ is concave upper semi-continuous on $M^*$, while $x\to b(t, v, x)$ is convex lower semi-continuous on $M$. Moreover, $b$ is locally Lipschitz continuous on $[0, \infty)\times M\times M^*$.

\item  The costs $b$, $c$ and $\tilde c$ are dual to each other in the following sense:
\begin{itemize}
\item For any $(v,x)\in M^*\times M$, we have $b(t,v, x)=\inf\{\langle v, y\rangle +c(t,y,x);\,  y\in M\}.$

\item For any $(y, x)\in M\times M$, we have $c(t,y,x)=\sup\{ b(t, v, x)-\langle v, y\rangle; v\in M^*\}.$

\item For any $(v,x)\in M^*\times M$, we have $b(t,v, x)=\sup\{\langle w, x\rangle - {\tilde c}(t, v, w); w\in M^*\}. $
\end{itemize}
 
\end{enumerate}
\end{prop}

The above  statements follow readily from properties of value functions that are consequences of the duality in Bolza's problem, which is a particular case of the following set up. Consider the path space 
$A^{2}_{M} = \{ u:[0,T] \rightarrow M; \,\dot{u} \in L^{2}_{M}  \}
$
equipped with the norm
\[
     \|u\|_{A^{^{2}}_{M}} = \left(\|u(0)\|_{M}^{2} +
     \int_0^T \|\dot{u}\|^{2} dt\right)^{\frac{1}{2}}.
     \]
      One way to represent the space  $A^{^{2}}_{M}$ is to identify it with the product space
$M \times L^{2}_{M}$, in such a way that its dual  $(A^{2}_{M})^*$ can also be
identified with $ M \times
L^{2}_{M}$  via the formula:
\begin{equation}\label{first.duality}
     \braket{u}{(p_{1},p_{0})}_{_{A^{2}_{M},M \times L_{M}^{2}}} =
      \braket{u(0)}{p_{1}} + \int_0^T \braket{\dot{u}(t)}{p_{0}(t)}dt,
\end{equation}
where $u\in A^{2}_{M}$ and $(p_{1},p_{0})\in M\times L_{M}^{2}$. We have the following duality formula.

\begin{prop} \label{Legendre.integral} Let $L$ be a time-dependent convex Lagrangian on $M\times M$ and let $\ell$ be a proper convex lower semi-continuous function on $M\times M$. Consider the Lagrangian  on $A^{2}_{M}\times (A^{2}_{M})^{*}=A^{2}_{M}\times (M \times L_{M}^{2})$ defined by  
 \begin{equation}
{\cal N} (u, p)= \int_0^TL(t, u(t)-p_0(t), -\dot u(t)) dt +\ell (u(0)-a, u(T))
\end{equation}
where  $u\in A^{2}_{M}$ and  $(p_0(t), a)\in L^2_M\times M$ represents an element $p$ in the dual of $A^{2}_{M}$. 
Then, for any $(v, q) \in A^{2}_{M}\times (A^{2}_{M})^{*}$ with $q$ of the form $(q_{0}(t), 0)$, we have 
 \begin{equation}
{\cal N}^* (q, v)= \int_0^TL^*(t, -\dot v(t), v(t)-q_0(t), ) dt +\ell^* (-v(0), v(T)).
\end{equation}
\end{prop}
 \noindent {\bf Proof:}  For $(v, q) \in A^{2}_{M}\times (A^{2}_{M})^{*}$ with $q$ represented by $(q_{0}(t), 0)$ write:
 \begin{eqnarray*}
  {\cal N}^{*}(q,v) &=& \sup_{p_{1} \in M} \sup_{p_{0} \in L_M^{2}} \sup_{u \in
     A^{2}_{M}} \left\{ \langle p_{1},v(0) \rangle + \int_0^T  \braket{p_{0}(t)}{\dot{v}(t)}+\braket{q_{0}(t)}{\dot{u}(t)}\, dt\right.\\
     &&    \quad \quad \quad \quad \left. 
     - \int_0^TL(t, u(t)-p_{0}(t),-\dot{u}(t))\, dt  -\ell (u(0)-p_{1},u(T)) \right\}.
\end{eqnarray*}
  Making a substitution $u(0) - p_{1} = a \in M$ and $u(t)-p_{0}(t)=y(t) \in L_M^2$,
we obtain
  \begin{eqnarray*}
    { {\cal N}}^{*}(q,v) &=& \sup_{a \in M} \sup_{y \in L_M^{2}} \sup_{u \in
     A^{2}_{M}} \left\{ \langle u(0)-a,v(0) \rangle -\ell (a,u(T))\right.\\
  && \quad \quad \left.+ \int_0^T \{ \braket{u(t)-y(t)}{\dot{v}(t)}+\braket{q_{0}(t)}{\dot{u}(t)}
     - L(t, y(t),-\dot{u}(t))\} dt \right\}
   \end{eqnarray*}
Since $\dot u$ and $\dot v\in L_M^2$, we have
\[
\int_0^T \braket{u}{\dot{v}}=- \int_0^T \braket{\dot u}{v} +
\braket{v(T)}{u(T)} - \braket{v(0)}{u(0)}
\] 
    which implies
 \begin{eqnarray*}
    { {\cal N}}^{*}(q,v) &=&\sup_{a \in M} \sup_{y \in L_M^{2}} \sup_{u \in
     A^{2}_{M}} \left\{ \langle -a,v(0) \rangle +\braket{v(T)}{u(T)}-\ell (a,u(T))\right.\\
  &&  \quad \quad \left. \int_0^T \left[ -\braket{y(t)}{\dot{v}(t)}+\braket{v(t)-q_{0}(t)}{-\dot{u}(t)} 
     - L(t, y(t),-\dot{u}(t))\right] dt \right\}.
  \end{eqnarray*}
Identify now $A_{M}^{^2}$ with $M\times L_M^2$
via the correspondence:
  \begin{eqnarray*}
     (b,r) \in M\times L_M^2 &\mapsto & b+\int_t^T r(s)\, ds\in
A_{M}^{^2}\\
     u\in A_{M}^{^2}   &\mapsto & \big( u(T),-\dot u(t)\big)\in M\times
L_M^2 .
\end{eqnarray*}
We finally obtain
  \begin{eqnarray*}
    { {\cal N}}^{*}(q,v) &=&\sup_{a \in M}\sup_{b \in M}  \{ \langle a,-v(0) \rangle +\braket{v(T)}{b}-\ell (a,b)\\
  &+& \sup_{y \in L_M^{2}} \sup_{r \in
     L^{2}_{M}}\int_0^T -\braket{y(t)}{\dot{v}(t)}+\braket{v(t)-q_{0}(t)}{r(t)} 
     - L(t, y(t),r(t)) dt\\
  &=&\int_{0}^{T}L^{*}(t, -\dot{v}(t),v(t)-q_{0}(t)) dt +\ell^{*}(-v(0), v(T)).
  \end{eqnarray*}
Now consider the minimization problems,
\begin{equation}\label{Bolza}
\hbox{$({\mathcal P})\qquad \inf\left\{
\int_0^TL(\gamma (s), {\dot \gamma}(s))\, ds +\ell (\gamma(0), \gamma (T)); \, \gamma \in  C^1([0, T), M)\right\},$}
\end{equation}
and  
\begin{equation}
\hbox{$({\tilde {\mathcal P}})\qquad \inf\left\{\int_0^T{\tilde L}(\gamma (s), {\dot \gamma}(s))\, ds +\ell^* (\gamma(0), -\gamma (T)); \, \gamma \in  C^1([0, T), M)\right\}$.}
\end{equation}
It is clear that 
$\inf({\mathcal P})\geq -\inf ({\tilde {\mathcal P}})$, and that we have equality under very mild conditions (see \cite{R1}). In this case, 
 there exists two paths $x(t), v(t)$ in $A^2_M$ satisfying 
 \begin{equation}
(\dot v(t), v(t))\in \partial L(x(t), \dot x(t)) \quad \hbox{for a.e. $t$,}
\end{equation} 
which can also be written in a dual form 
\begin{equation}
({\dot x}(t), x(t))\in \partial {\tilde L}(v(t), \dot v(t)) \quad \hbox{for a.e. $t$,}
\end{equation} 
or in a Hamiltonian form as
\begin{eqnarray}
\dot x(t)&\in& \partial_v H(x(t), v(t))\\
-\dot v(t)&\in& {\tilde \partial}_x H(x(t), v(t)), 
\end{eqnarray} 
coupled with the boundary conditions
\begin{equation}
(v(0), -v(T))\in \partial \ell (x(0), x(T)).
\end{equation}
 Note that $c_T$, ${\tilde c}_T$, $b_T$ and all value functions $V_g$ for $g$ convex, can be written as in (\ref{Bolza}), and therefore inherit whatever the duality $\inf({\mathcal P})=-\inf ({\tilde {\mathcal P}})$ provides.  

\begin{prop} The value function $V_g$ is expressed in terms of the costs $c$ and $b$ by the following formulae:
\begin{enumerate}
\item $V_g(t,x)=\inf\{g(y)+c(t, y, x);\, y\in M\}.$

\item If $g$ is convex and lower semi-continuous, then $V_g(t,x)=\sup\{b(t, v, x)-g^*(v); \, v\in M^*\}.$

\item For each $t\in [0, +\infty)$, the Legendre transform of the convex function $x\to V_g(t, x)$ on $M$ is the convex function  $w\to \tilde V_{g^*}(t, w)$  given by
\begin{equation*}
\left\{ \begin{array}{lll}
{\tilde V}_{g^*}(t,w)&=&\inf\{g^*(\gamma (0))+\int_0^t{\tilde L}(s, \gamma (s), {\dot \gamma}(s))\, ds; u\in C^1([0, T), M);   \gamma(t)=w\},\\
{\tilde V}_{g^*}(0,w)&=&g^*(w).
\end{array}\right.
\end{equation*}

\item For each $t$, the graph of the subgradient $\partial V_g(t, \cdot)$, i..e., 
$$\Gamma_g(t)=\{(x, v); v\in \partial V_g(t, x)\}$$
is a globally Lipchitzian manifold of dimension $n$ in $M \times M^*$, which depends continuously on $t$.

\end{enumerate}
\end{prop}
Note that $b(t, v, x)=V_g(t,x)$, when $g_v(y)=\langle v, y\rangle$. In this case, $g_v^*(u)=0$ if $u=v$  and $+\infty$ if $u\neq v$, which yields that the Legendre dual of $x\to V_g(t, x)=b(t, v, x)$ is $w\to {\tilde c}(t, v , w)$. One can also deduce the following.

\begin{prop} \label{key} Let $g$ be a proper convex lower semi-continuous function on $M$. 
\begin{enumerate}
\item If a Hamiltonian trajectory $(x(t), v(t))$ over $[0, T ]$ starts with $v(0) \in \partial g(x(0))$, then $v(t) \in \partial V_g (t, x(t))$ for all $t \in [0, T]$.

\item Moreover, a pair of arcs $x(t)$ and $v(t)$ gives a
Hamiltonian trajectory over $[0, T]$ that starts in the graph of $\partial g$ and ends at a point $x, w$ in the graph of $\partial V_g(T,\cdot )$ if and only if
$x(t)$ is optimal in the minimization problem  that defines $V_g(t, x)$ and
$v(t)$ is optimal in the minimization problem that defines $\tilde V_g(t, w)$.

\item In particular, the following properties are equivalent:
\begin{enumerate}
\item $(-v,w)\in \partial_{y,x}c_T(y, x)$;

\item $w \in \partial_{x}b_T (v, x)$ and $y \in {\tilde \partial}_{v}b_T(v, x)$.

\item There is a Hamiltonian trajectory $(\gamma (t), \eta (t))$ over $[0, T]$ starting at $(y, v)$ and ending at $(x, w)$.

\end{enumerate}
\end{enumerate}
\end{prop}
We finally recall the twist condition.

\begin{defn} A cost function $c$ satisfies {\it the twist condition} if for each $y\in M$, we have $x=x'$ whenever the differentials  $\partial_yc(y,x)$ and $\partial_yc(y,x')$ exist and are equal.
\end{defn}
In view of the above proposition, $c_T$ satisfies the twist condition if there is at most one Hamiltonian trajectory starting at a given initial state $(v, y)$, while the cost $b_T$ satisfies the twist condition if for any given states $(v, w)$, there is at most one Hamiltonian trajectory starting at $v$ and ending at $w$.  

\section{The Hopf-Lax formulas on Wasserstein space}

We shall frequently use the following results of Brenier \cite{B1} that we summarize in the following proposition.

\begin{prop} (Brenier) If $\nu_0$ (resp., $\mu_0$) is a probability measure on $M$ (resp., $M^*$) such that $\nu_0$ is absolutely continuous measure with respect to Lebesgue measure, then there exists a unique (up to a constant) convex function $h$ (resp., concave function $g$) on $M$ such that $(\nabla h)_\#\nu_0=\mu_0$ (resp., $(\nabla g)_\#\nu_0=\mu_0$) and 
$$
\hbox{${\overline W}(\nu_0, \mu_0)=\int_M\langle \nabla h(x), x\rangle \, d\nu_0(x)$,\quad  resp., ${\underline W}(\nu_0, \mu_0)=\int_M\langle \nabla g(x), x\rangle \, d\nu_0(x)$.}
$$
Moreover, ${\underline W}(\nu_0, \mu_0)=-{\overline W}({\tilde \nu_0}, \mu_0),$ where 
${\tilde \nu_0}(A)=\nu_0(-A)$.
\end{prop} 
Note that if the convex function $h$ solves ${\overline W}({\tilde \nu_0}, \mu_0)$, then the concave function $g(x)=-h(-x)$ solves ${\underline W}(\nu_0, \mu_0)$ and vice-versa.\\

We shall need the following simple lemma regarding Kantorovich potentials.

\begin{lem} \label{pots} Let $g$ be a function on $M^*$, and let  $h, k$ be functions on $M$. 
\begin{enumerate}
\item If $h(x)-g(v) \geq b_T(v, x)$ on $M^*\times M$, then $h^*(w)+g(v) \leq {\tilde c}_T(v, w)$ on $M^*\times M^*$. 
\item If $h(x)-g(v) \leq b_T(v, x)$ on $M^*\times M$, then $h(x)+g^*(-y) \leq c(y, x)$ on $M\times M$.
\item If $h(x)-g(y) \leq c(y, x)$ on $M\times M$, then $h(x)-(-g)^*(-v)\leq b_T(v, x)$ on $M^*\times M$.
\end{enumerate}
\end{lem}

\noindent {\bf Proof:} 1) Since for any $(v, w)\in M^*\times M^*$, we have 
$${\tilde c}(t,v,w)=\sup\{\langle w, x\rangle - b(t, v, x); x\in M\}. $$
It follows that for any $y\in M$, 
\begin{eqnarray*}
{\tilde c}_T(v,w)\geq  \langle w, y\rangle - b(t, v, y)\geq \langle w, y\rangle +g(v)-h(y) 
\end{eqnarray*}
hence
\[
h^*(w) +g(v) \leq {\tilde c}_T(v, w).  
\]
2)  Since for any $(y, x)\in M\times M$, we have $c(t,y,x)=\sup\{ b(t, v, x)-\langle v, y\rangle; v\in M^*\}$,  it follows that for any $y\in M$, we have $-g(v)+h(x)-\langle v, y\rangle \leq c(x, y)$, 
that is 
\[
h(x) +g^*(-y) \leq c(y,x).
\]
3) Since $b(t,v, x)=\inf\{\langle v, y\rangle +c(t,y,x);\,  y\in M\}$, it follows that for any $(v, x)\in M^*\times M$, and any $\epsilon>0$, there exists $y_\epsilon \in M$ such that 
\[
c_T(y_0, x) +\langle v, y_0\rangle \leq b_T(v, x) +\epsilon,
\]
hence $h(x)-g(y_0)+\langle v, y_0\rangle \leq b_T(v, x) +\epsilon$, which means that 
\[
h(x) +\inf\{-g(y)+\langle v, y\rangle\} \leq b_T(v,x),
\]
that is $h(x)-(-g)^*(-v)\leq b_T(v, x)$.\\

\noindent {\bf Proof of Theorem \ref{interpol}:} First, we prove (\ref{HL.one}), that is if $\mu_0$ is absolutely continuous with respect to Lebesgue measure, then   
\begin{equation}\label{two.again}
{\underline B}_T(\mu_0,\nu_T)=\inf\{ C_T(\nu, \nu_T) +{\underline W}(\mu_0, \nu);\, \nu\in {\mathcal P}(M)\}.  
\end{equation}
We note that for any probability measure $\nu$ on $M$, we have 
\begin{equation}\label{inequality.again}
{\underline B}_T(\mu_0,\nu_T)\leq  C_T(\nu, \nu_T) +{\underline W}(\mu_0, \nu).  
\end{equation}
Indeed, since $\mu_0$ is assumed to be absolutely continuous with respect to Lebesgue measure, Brenier's theorem yields a concave function $k$ on $M^*$ that is differentiable $\mu_0$-almost everywhere such that  $(\nabla k)_\# \mu_0=\nu$, and  
$${\underline W}(\mu_0, \nu)=\int_{M^*}\langle v, \nabla k(v) \rangle \, d\mu_0(v).
$$
 Also let $\pi_0$ be an optimal transport plan for 
$C_T(\nu, \nu_T)$, that is $\pi_0\in {\mathcal K}(\nu, \nu_T)$ such that
\[
C_T(\nu, \nu_T)=\int_{M\times M} c_T(y,x)\, d\pi_0(y,x).
\]
Let  ${\tilde \pi}_0:= S_\#\pi_0$, where $S(y, x)=(\nabla {\tilde k}(y), x)$, where ${\tilde k}$ is the concave Legendre transform of $k$. It is a transport plan in ${\mathcal K}( \mu_0, \nu_T)$. Since $ b_T(v, x)\leq  c_T(\nabla k(v), x )+\langle \nabla k(v), v\rangle$ for every $v\in M^*$, we have 
\begin{eqnarray*}
{\underline B}_T(\mu_0,\nu_T)&\leq&\int_{M^*\times M} b_T(v,x)\, d{\tilde\pi}_0(v, x)\\
&\leq& \int_{M^*\times M}\{c_T(\nabla k(v), x )+\langle \nabla k(v), v\rangle \}d{\tilde \pi}_0(v, x)\\
&= &\int_{M\times M}c_T(y,x) d\pi_0(y,x)+\int_{M}\langle \nabla k(v), v\rangle \, d \mu_0(v)\\
&=& C_T(\nu, \nu_T)+{\underline W}(\mu_0, \nu). 
\end{eqnarray*}
To prove the reverse inequality, use standard Monge-Kantorovich theory to write
\begin{eqnarray*}
{\underline B}_T(\mu_0,\nu_T)&=&\inf\{\int_{M^*\times M} b_T(v, x)\, d\pi(v, x);\, \pi \in {\mathcal K}(\mu_0, \nu_T)\}\\
&=&\sup\{\int_{M} h(x)\, d\nu_T(x)-\int_{M^*}g(v)\, d\mu_0(v);\, h(x)-g(v) \leq b_T(v, x)\}.
\end{eqnarray*} 
Since the supremum can be taken over all admissible Kantorovich pairs $(g, h)$ of functions, i.e. those satisfying
the relations
$$
g(v)=\sup_{x\in M} h(x)-b_T(v, x)
\text{\quad  and \quad } 
h(x)=\inf _{v\in M^*} b_T(v, x))+g(v), 
$$
we can assume that the initial Kantorovic potential $g$ is convex. Since the cost function $b_T$ is continuous, the infimum ${\underline B}_T(\mu_0,\nu_T)$ is attained at some probability measure $\pi_0\in {\mathcal K}(\mu_0, \nu_T)$. Moreover, the supremum in the dual problem  is attained at some pair $(g, h)$ of admissible Kantorovich functions. It follows that $\pi_0$ is supported on the set
\[
{\mathcal O}:=\{(v, x)\in M^*\times M; \, b_T(v, x)=h(x)-g(v)\}.
\]
We now exploit the convexity of $g$, and use the fact that 
for each $(v, x) \in {\mathcal O}$, the function $w\to -g(w)+h(x)-b_T(w,x)$ attains its maximum at $v$, which means that 
\[
-\nabla g(v)\in {\tilde \partial}_v b_T(v, x), 
\]
where for a concave function $f$, we write ${\tilde \partial} f:=-\partial (-f)$, the latter being negative the subdifferential of the convex function $-f$.
But since $y\to -c_T(y, x)$ is the concave Legendre transform of $v\to b_T(v, x)$ with respect to the $v$-variable, we then have 
\begin{equation}
b_T(v, x)-c_T(-\nabla g(v), x)=\langle v, -\nabla g(v)\rangle.
\end{equation}
Integrating with $\pi_0$, we get since $\pi_0\in {\mathcal K}(\mu_0, \nu_T)$,
\begin{equation}
\int_{M^*\times M}b_T(v, x)\, d\pi_0-\int_{M^*\times M}c_T(-\nabla g(v), x)d\pi_0=\int_{M^*}\langle v, -\nabla g(v)\rangle\, d\mu_0.
\end{equation}
Letting $\nu_0=(-\nabla g)_\#\mu_0$, we obtain that 
\begin{equation}\label{crucial.two}
{\underline B}_T(\mu_0,\nu_T)-\int_{M^*\times M}c_T(-\nabla g(v), x)d\pi_0={\underline W}(\mu_0, \nu_0).
 \end{equation}
We now prove that 
\begin{equation}\label{magic.two}
\int_{M^*\times M}c_T(-\nabla g(v), x)d\pi_0=C_T(\nu_0, \nu_T).
\end{equation}
Indeed, we have 
\begin{eqnarray*}
\int_{M^*\times M}c_T(-\nabla g(v), x)d\pi_0\geq C_T(\nu_0, \nu_T), 
\end{eqnarray*}
since 
the measure $\pi=S_\#\pi_0$, where $S(v,x)=(-\nabla g(v), x)$ has marginals $\nu_0$ and $\nu_T$ respectively. On the other hand, 
 (\ref{crucial.two}) yields 
\begin{eqnarray*}
\int_{M^*\times M}c_T(-\nabla g(v), x)d\pi_0&=&{\underline B}_T(\mu_0,\nu_T)+\int_{M^*}\langle v, \nabla g(v)\rangle\, d\mu_0(v)\\
&=&\int_Mh(x)\, d\nu_T(x)-\int_{M^*}g(v)\, d\mu_0(v)+\int_{M^*} g^*(\nabla g(v))\, d \mu_0(v)+ \int_{M^*}g(v)\, d\mu_0(v)\\
&=&\int_Mh(x)\, d\nu_T(x)+\int_{M^*}g^*(-y) d\nu_0(y).
\end{eqnarray*}
Moreover, since $h(x)-g(v) \leq b(v,x)$ for all $(x, v)\in M\times M^*$, it follows from Lemma \ref{pots} 
that 
\begin{equation}
h(x) +g^*(-y) \leq c_T(y,x), 
\end{equation}
which means that the couple $(-g^*(-y), h(x))$ is an admissible Kantorovich pair for the cost $c_T$. Hence,
\begin{eqnarray*}
C_T(\mu_0, \mu_T)&\leq& \int_{M^*\times M}c_T(-\nabla g(v), x)d\pi_0\\
&=&\int_Mh(x)\, d\nu_T(x)+\int_{M^*}g^*(-y) d\nu_0(y)\\
&\leq&\sup\{\int_{M} \phi_1(x)\, d\nu_T(x)-\int_{M}\phi_0(y)\, d\nu_0(y);\, \phi_1(x)- \phi_0(y) \leq  c_T(y, x)\}\\
&=&C_T(\mu_0, \mu_T).
\end{eqnarray*}
It follows that
\begin{equation}
{\underline B}_T(\mu_0,\nu_T)=\inf\{ {\underline W}(\mu_0, \nu)+C_T(\nu, \nu_T);\, \nu\in {\mathcal P}(M)\}.  
\end{equation}
and the infimum is attained by the measure $\nu_0$. Note that the final optimal Kantorovich potential for $C_T(\nu_0, \nu_T)$ is $y\to -g^*(-y)$, hence is concave.

Similarly, we can show (\ref{HL.two}). First, for any probability measure $\mu$ on $M^*$, we have 
\begin{equation}\label{inequality}
{\overline B}_T(\mu_0,\nu_T)\geq  {\overline W}(\nu_T, \mu)-{\tilde C}_T(\mu_0, \mu).  
\end{equation}
Indeed, since $\nu_T$ is assumed to be absolutely continuous with respect to Lebesgue measure, Brenier's theorem yields a convex function $h$ that is differentiable $\mu_T$-almost everywhere on $M$ such that $(\nabla h)_\#\nu_T=\mu$, and  ${\overline W}(\nu_T, \mu)=\int_{M}\langle x, \nabla h(x)\rangle \, d\nu_T(x)$. Let $\pi_0$ be an optimal transport plan for 
${\tilde C}_T(\mu_0, \mu)$, that is $\pi_0\in {\mathcal K}(\mu_0, \mu)$ such that
\[
{\tilde C}_T(\mu_0, \mu)=\int_{M^*\times M^*} {\tilde c}_T(v,w)\, d\pi_0(v, w).
\]
Let ${\tilde \pi}_0:= S_\#\pi_0$, where $S(v,w)=(v, \nabla h^*(w))$, which is a transport plan in ${\mathcal K}(\mu_0, \nu_T)$. Since $ b_T(v, y)\geq \langle \nabla h(x), y\rangle -{\tilde c}_T(v, \nabla h(x) )$ for every $(y, x, v)\in M\times M\times M^*$, we have 
\begin{eqnarray*}
{\overline B}_T(\mu_0,\nu_T)&\geq&\int_{M^*\times M} b_T(v,x)\, d{\tilde \pi}_0(v, x)\\
&\geq& \int_{M^*\times M}\{ \langle \nabla h(x), x\rangle -{\tilde c}_T(v, \nabla h(x) )\} d{\tilde \pi}_0(v, x)\\
&= &\int_{M}\langle x, \nabla h(x)\rangle \, d\nu_T(x) -\int_{M^*\times M^*}{\tilde c}_T(v, w)\} d\pi_0(v, w)\\
&=& {\overline W}(\nu_T, \mu)-{\tilde C}_T(\mu_0, \mu).  
\end{eqnarray*}
To prove the reverse inequality, we use standard Monge-Kantorovich theory to write
\begin{eqnarray*}
{\overline B}_T(\mu_0,\nu_T)&=&\sup\big\{\int_{M^*\times M} b_T(v, x)\, d\pi(v, x);\, \pi \in {\mathcal K}(\mu_0, \nu_T)\big\}\\
&=&\inf\big\{\int_{M} h(x)\, d\nu_T(x)-\int_{M^*}g(v)\, d\mu_0(v);\, h(x)-g(v) \geq b_T(v, x)\big\},
\end{eqnarray*} 
where the infimum is taken over all admissible Kantorovich pairs $(g, h)$ of functions, i.e. those satisfying
the relations
$$
g(v)=\inf_{x\in M}h(x)- b_T(v, x)
\text{\quad  and \quad } 
h(x)=\sup _{v\in M^*} b_T(v, x))+g(v)
$$
Note that $h$ is convex.
Since the cost function $b_T$ is continuous, the supremum ${\overline B}_T(\mu_0,\nu_T)$ is attained at some probability measure $\pi_0\in {\mathcal K}(\mu_0, \nu_T)$. Moreover, the infimum in the dual problem  is attained at some pair $(g, h)$ of admissible Kantorovich functions. It follows that $\pi_0$ is supported on the set
\[
{\mathcal O}:=\{(v, x)\in M^*\times M; \, b_T(v, x)=h(x)-g(v)\}.
\]
We now exploit the convexity of $h$, and use the fact that 
for each $(v, x) \in {\mathcal O}$, the function $y\to h(y)-g(v)-b_T(v,y)$ attains its minimum at $x$, which means that 
\[
\nabla h(x)\in {\partial}_x b_T(v, x).
\]
But since $\tilde c_T$ is the Legendre transform of $b_T$ with respect to the $x$-variable, we then have 
\begin{equation}
b_T(v, x)+{\tilde c}_T(v, \nabla h(x))=\langle x, \nabla h(x)\rangle \,\, \hbox{on ${\mathcal O}$}.
\end{equation}
Integrating with $\pi_0$, we get since $\pi_0\in {\mathcal K}(\mu_0, \nu_T)$,
\begin{equation}
\int_{M^*\times M}b_T(v, x)\, d\pi_0+\int_{M^*\times M}{\tilde c}_T(v, \nabla h(x))d\pi_0=\int_{M}\langle x, \nabla h(x)\rangle\, d\nu_T.
\end{equation}
Letting $\mu_T=\nabla h_{\#}\nu_T$, we obtain that 
\begin{equation}\label{crucial}
{\overline B}_T(\mu_0,\nu_T)+\int_{M^*\times M}{\tilde c}_T(v, \nabla h(x))d\pi_0={\overline W}(\nu_T, \mu_T), \end{equation}
where 
\[
{\overline W}(\nu_T, \mu_T)=\sup\{\int_{M\times M^*} \langle x,v\rangle\, d\pi;\, \pi \in {\mathcal K}(\nu_T, \mu_T)\}.
\]
Note that we have used here that $h$ is convex to deduce that ${\overline W}(\nu_T, \mu_T)=\int_{M}\langle x, \nabla h(x)\, d\mu_T$ by the uniqueness in Brenier's decomposition. We now prove that 
\begin{equation}\label{magic}
\int_{M^*\times M}{\tilde c}_T(v, \nabla h(x))d\pi_0={\tilde C}_T(\mu_0, \mu_T), 
\end{equation}
where 
\[
{\tilde C}_T(\mu_0, \mu_T):=\inf\{\int_{M^*\times M^*} {\tilde c}_T(v, w)d\pi ;\, \pi \in {\mathcal K}(\mu_0, \mu_T)\}.
\]
Indeed, we have 
\begin{eqnarray*}
\int_{M^*\times M}{\tilde c}_T(v, \nabla h(x))d\pi_0&\geq& {\tilde C}_T(\mu_0, \mu_T),
\end{eqnarray*}
since the measure $\pi=S_\#\pi_0$, where $S(v,x)=(v, \nabla h(x))$ has marginals $\mu_0$ and $\mu_T$ respectively. On the other hand, 
 (\ref{crucial}) yields 
\begin{eqnarray*}
\int_{M^*\times M}{\tilde c}_T(v, \nabla h(x))d\pi_0&=&\int_{M}\langle x, \nabla h(x)\rangle\, d\nu_T(x)-\int_{M^*\times M}b_T(v, x)\, d\pi_0\\
&=&\int_{M}h^*(\nabla h(x)) d\nu_T(x)+ \int_Mh(x)\, d\nu_T(x)+\int_{M^*}g(v)\, d\mu_0(v)-\int_Mh(x)\, d\nu_T(x)\\
&=&\int_{M^*}h^*(w) d\mu_T(w)+\int_{M^*}g(v)\, d\mu_0(v).
\end{eqnarray*}
Moreover, since $h(x)-g(v) \geq b(v, x)$, we have by Lemma \ref{pots}, that 
$$
h^*(w) +g(v) \leq {\tilde c}_T(v, w),  
$$
which means that the couple $(-g, h^*)$ is an admissible Kantorovich pair for the cost ${\tilde c}_T$. Hence,
\begin{eqnarray*}
{\tilde C}_T(\mu_0, \mu_T)&\leq& \int_{M^*\times M}{\tilde c}_T(v, \nabla h(x))d\pi_0\\
&=&\int_{M}h^*(w) d\mu_T(w)+\int_{M^*}g(v)\, d\mu_0(v)\\
&\leq&\sup\{\int_{M^*} \phi_T(w)\, d\mu_T(w)-\int_{M^*}\phi_0(v)\, d\mu_0(v) ;\,  \phi_T(w)-\phi_0(v)  \leq {\tilde c}_T(v, w)\}\\
&=&{\tilde C}_T(\mu_0, \mu_T).
\end{eqnarray*}
It follows that
${\overline B}_T(\mu_0,\nu_T)={\overline W}(\nu_T, \mu_T)-{\tilde C}_T(\mu_0, \mu_T)$. 
In other words, the supremum in (\ref{inequality}) is attained by the measure $\mu_T$. Note that the final optimal Kantorovich potential for ${\tilde C}_T(\mu_0, \mu_T)$ is $h^*$, hence is convex.\\

\section{The reverse interpolation}

While the cost $c_T$ is itself jointly convex in both variables, one cannot deduce much in terms of the convexity or concavity of the corresponding Kantorovich potentials. 
We investigate here what happens in such situations. Note that we get a concave initial Kantorovich potential for ${C}_T(\nu_0, \nu_T)$, when $\nu_0$ is  obtained via the factorization of a ballistic optimal transport problem. Actually, the following somewhat converse statement holds true. 

\begin{thm}\label{endpoints} Assume $M=\R^d$ and that $L$ satisfies hypothesis (A1), (A2) and (A3). Assume $\nu_0$ and $\nu_T$ are probability measures on $M$ such that $\nu_0$ is absolutely continuous with respect to Lebesgue measure. Then, the initial Kantorovich potential of $C_T(\nu_0, \nu_T)$ is concave if and only if 
\begin{equation}\label{three}
C_T(\nu_0, \nu_T)=\sup\{{\underline B}_T(\mu,\nu_T)-{\underline W}(\nu_0, \mu);\, \mu\in {\mathcal P}(M^*)\}.  
\end{equation}
\end{thm}
\noindent {\bf Proof:} Again, we show first that 
\begin{equation}\label{three.prime}
C_T(\nu_0, \nu_T)\geq \sup\{{\underline B}_T(\mu,\nu_T)- {\underline W}(\nu_0, \mu);\, \mu\in {\mathcal P}(M^*)\}.  
\end{equation}
Indeed, if $\nu$ is any probability measure on $M$, use
again Brenier's theorem to find a concave function $h$ on $M$ such that  $(\nabla h)_\#\nu_0 =\mu$ and ${\underline W}(\nu_0, \mu)=\int_{M}\langle x, \nabla h(x)\rangle \, d\nu_0(x)$. Also let $\pi_0$ be an optimal transport plan for $C_T(\nu_0, \nu_T)$. 
Since $ c(t, y,x)\geq  b(t, \nabla h(y), x)- \langle \nabla h(y), y\rangle$ for every $x,y\in M$, and since 
${\tilde \pi_0}:=S_\#\pi_0$ where $S( y, x)= (\nabla h(y), x)$, is in ${\mathcal K}(\mu, \nu_T)$, we have  
\begin{eqnarray*}
C_T(\nu_0, \nu_T)&=&\int_{M\times M} c_T(y, x)\, d\pi_0\\
&\geq&
 \int_{M\times M}\{b_T(\nabla h(y),x)-\langle \nabla h(y), y\rangle 
 \} d\pi_0(y, x)\\
 &=& \int_{M\times M}b_T(\nabla h(y),x) d \pi_0(y, x)-{\underline W}(\nu_0, \mu)\\
 &=&\int_{M^*\times M}b_T(v,x) d{\tilde \pi}_0(v, x)-{\underline W}(\nu_0, \mu)\\
  &\geq &{\underline B}_T(\mu,\nu_T)- {\underline W}(\nu_0, \mu).
 \end{eqnarray*}
To prove Theorem \ref{endpoints}, we start by proving the reverse inequality and attainment in the case where the initial Kantorovich potential is concave.  Indeed, write 
\begin{eqnarray*}
C_T(\nu_0,\nu_T)&=&\inf\{\int_{M\times M} c(y, x)\, d\pi(y, x);\, \pi \in {\mathcal K}(\nu_0, \nu_T)\}\\
&=&\sup\{\int_{M} h(x)\, d\nu_T(x)-\int_{M}g(y)\, d\nu_0(y);\, h(x)-g(y) \leq c_T(y, x)\}.
\end{eqnarray*} 
Since the cost function $c_T$ is continuous, the infimum $C_T(\nu_0,\nu_T)$ is attained at some probability measure $\pi_0\in {\mathcal K}(\nu_0, \nu_T)$. Moreover, the infimum in the dual problem  is attained at some pair $(g, h)$ of admissible Kantorovich functions.  It follows that $\pi_0$ is supported on the set
\[
{\mathcal O}:=\{(y, x)\in M \times M; \, c_T(y, x)=h(x)-g(y)\}
\]
Assuming $g$ concave, use the fact that  for each $(y, x) \in {\mathcal O}$, the function $z\to h(x)-g(z) -c_T(z,x)$ attains its maxmum at $y$, to deduce  that 
\[
-\nabla g(y)\in {\partial}_y c_T(y, x).
\]
Assuming $g$ concave and since $b(t,v, x)=\inf\{\langle v, z\rangle +c(t,z,x);\,  z\in M\}$, this means that for $(y, x) \in {\mathcal O}$,
\begin{equation}
c_T(y, x)=b_T(\nabla g(y), x)-\langle \nabla g(y), y\rangle.
\end{equation}
Integrating with $\pi_0$, we get since $\pi_0\in {\mathcal K}(\nu_0, \nu_T)$,
\begin{equation}
\int_{M\times M}c_T(y, x)\, d\pi_0=\int_{M\times M}b_T(\nabla g(y), x)\, d\pi_0-\int_{M}\langle \nabla g(y), y\rangle\, d\nu_0.
\end{equation}
Letting $\mu_0=(\nabla g)_{\#}\nu_0$, and since $g$ is concave, we obtain that 
\begin{equation}\label{crucial.2}
C_T(\nu_0,\nu_T)=\int_{M\times M}b_T(\nabla g(y), x)\, d\pi_0 -{\underline W}(\nu_0, \mu_0).
 \end{equation}
We now prove that 
\begin{equation}\label{magic.2}
\int_{M\times M}b_T(\nabla g(y), x)\, d\pi_0(y,x)={\underline B}_T(\mu_0, \nu_T). 
\end{equation}
Indeed, we have 
\begin{eqnarray*}
\int_{M\times M}b_T(\nabla g(y), x)\, d\pi_0&\geq& {\underline B}_T(\mu_0, \nu_T),
\end{eqnarray*}
since the measure $\pi=S_\#\pi_0$ where $S(y, x)=(\nabla g(y), x)$ has $\mu_0$ and $\nu_T$ as marginals.
On the other hand, 
(\ref{crucial.2}) yields 
\begin{eqnarray*}
\int_{M\times M}b_T(\nabla g(y), x)\, d\pi_0&=&\int_{M\times M}c_T(y, x)\, d\pi_0+\int_{M}\langle y, \nabla g(y)\rangle\, d\nu_0(y)\\
&=&\int_Mh(x)\, d\nu_T(x)-\int_{M}g(y)\, d\nu_0(y)-\int_{M}(-g)^*(-\nabla g(y)) d\nu_0(y)+ \int_Mg(y)\, d\nu_0(y)\\
&=&\int_{M}h(x)\, d\nu_T(x)-\int_{M^*}(-g)^*(-v) d\mu_0(v).
\end{eqnarray*}
Moreover, by Lemma \ref{pots}, we have that
$h(x)-(-g)^*(-v) \leq b_T(v,x)$, that is the couple $((-g)^*(-v), h (x))$ is an admissible Kantorovich pair for the cost $b_T$. It follows that
\begin{eqnarray*}
{\underline B}_T(\mu_0, \nu_T)&\leq&\int_{M\times M}b_T(\nabla g(y), x)\, d\pi_0\\
&=&\int_{M}h(x)\, d\nu_T(x)-\int_{M}(-g)^*(-v) d\mu_0(v)\\
&\leq&\sup\{\int_{M} \phi_T(x)\, d\mu_T(x)-\int_{M^*}\phi_0(v)\, d\mu_0(v);\,  \phi_T(x)-\phi_0(v)  \leq b_T(v, x)\}\\
&=&{\underline B}_T(\mu_0, \nu_T),
\end{eqnarray*}
and 
$C_T(\nu_0, \nu_T)={\underline B}_T(\mu_0,\nu_T)-{\underline W}(\nu_0,  \mu_0). $
In other words, the supremum in (\ref{three}) is attained by the measure $\mu_0$.

\begin{cor} Assume $M=\R^d$ and that $L$ satisfies hypothesis (A1), (A2) and (A3). Assume $\nu_0$ and $\nu_T$ are probability measures on $M$ such that $\nu_0$ is absolutely continuous with respect to Lebesgue measure, and that the initial Kantorovich potential of $C_T(\nu_0, \nu_T)$ is concave. If $b_T$ satisfies the twist condition, then there exists a map $X_0^T: M^*\to M$ and a concave function $g$ on $M$ such that 
\begin{equation}
C_T(\nu_0, \nu_T)=\int _{M} c_T (y, X_0^T\circ \nabla g(y)) d\nu_0(y).
\end{equation}
\end{cor}
\noindent{\bf Proof:} In this case,
\begin{equation}\label{three.bis}
C_T(\nu_0, \nu_T)={\
\underline B}_T(\mu_0,\nu_T)-{\underline W}(\nu_0, \mu_0), 
\end{equation}
for some probability measure $\mu_0$ on $M^*$. Let $g$ be the concave function on $M$ such that 
$(\nabla g)_\#\nu_0=\mu_0$ and 
\[
{\underline W}(\nu_0, \mu_0)=\int_M \langle \nabla g(y), y\rangle d\nu_0(y).
\]
Since $b_T$ satisfies the twist condition, there exists a map $ X_0^T: M^*\to M$  such that $(X_0^T)_\#\mu_0=\nu_T$ and 
\begin{equation}
{\underline B}_T(\mu_0, \nu_T)=\int _{M^*} b_T (v, X_0^Tv) d\mu_0(v).
\end{equation}
Note that the infimum $C_T(\nu_0, \nu_T)$  is attained at some probability measure $\pi_0\in {\mathcal K}(\nu_0, \nu_T)$ and that $\pi_0$ is supported on a subset ${\mathcal O}$ of $M\times M$ such that 
for $(y, x) \in {\mathcal O}$,
\begin{equation*}
c_T(y, x)=b_T(\nabla g(y), x)-\langle \nabla g(y), y\rangle.
\end{equation*}
Moreover, $C_T(\nu_0,\nu_T)=\int_{M\times M}b_T(\nabla g(y), x)\, d\pi_0 -{\underline W}(\nu_0, \mu_0)$, and 
\[
\int_{M\times M}b_T(\nabla g(y), x)\, d\pi_0={\underline B}_T(\mu_0, \nu_T)=\int _{M^*} b_T (v, X_0^Tv) d\mu_0(v)=\int _{M} b_T (\nabla g(y), X_0^T\circ \nabla g(y)) d\nu_0(y).
\]
Since $b_T$ satisfies the twist condition, it follows that for any $(y, x)\in {\mathcal O}$, we have that 
$x=X_0^T\circ \nabla g(y)$ from which follows that  $C_T(\nu_0, \nu_T)=\int _{M} c_T (y, X_0^T\circ \nabla g(y)) d\nu_0(y).$\\

\noindent {\bf Proof of Corollary \ref{cor1}:} The cost $c(x-y)$ corresponds to $c_1(y, x)$, where the Lagrangian is $L(x, v)=c(v)$, that is 
\begin{equation}\label{v}
c_1(y,x)=\inf\{\int_0^1c({\dot \gamma}(t))\, dt; \gamma\in C^1([0, 1), M);  \gamma(0)=y, \gamma(1)=x\}=c(x-y).
\end{equation}
It follows from (\ref{three}) that there is a probability measure $\mu_0$ on $M^*$ such that 
$C_1(\nu_0, \nu_1)={\underline B}_1(\mu_0,\nu_1)-{\underline W}(\nu_0, \mu_0)$. But
in this case,
$b_1(v, x)=\inf\{\langle v, y\rangle +c(x-y);\,  y\in M\}=\langle v, x\rangle-c^*(v)$, hence 
\begin{equation}\label{ng}
C_1(\nu_0, \nu_1)={\underline B}_1(\mu_0,\nu_1)-{\underline W}(\nu_0, \mu_0)={\underline W}(\mu_0, \nu_1)-\int_{M^*}c^*(v)\, d\mu_0(v)-{\underline W}(\nu_0, \mu_0). 
\end{equation}
In other words, 
\begin{equation}\label{ng}
C_1(\nu_0, \nu_1)+K={\underline W}_1(\mu_0,\nu_1).-{\underline W}(\nu_0, \mu_0),
\end{equation}
where $K$ is the constant $\int_{M^*}c^*(v)\, d\mu_0(v)$. 

Apply Brenier's theorem \cite{B1} twice to find concave functions $\phi_0: M\to \R$ and $\phi_1: M^*\to \R$ such that $(\nabla \phi_0)_\#\nu_0=\mu_0$, $(\nabla \phi_1)_\#\mu_0=\nu_1$ and
\[
\hbox{${\underline W}(\nu_0, \mu_0)=\int_{M}\langle y, \nabla \phi_0(y)\rangle \, d \nu_0(y)$
\quad and \quad 
$
{\underline W}( \mu_0, \nu_1)=\int_{M^*}\langle v, \nabla \phi_1(v)\rangle \, d\mu_0(v).$}
\]
It follows from the preceeding corollary that 
\begin{eqnarray*}
C_1(\nu_0, \nu_1)+K=\int _{M} c_1 (y, \nabla \phi_1\circ \nabla \phi_0(y)) d\nu_0(y)=\int _{M} c ( \nabla \phi_1\circ \nabla \phi_0(y)-y) d\nu_0(y).
\end{eqnarray*}
Note also that 
\begin{eqnarray*}
C_1(\nu_0, \nu_1)+K&=&\int_{M}\langle v, \nabla \phi_1(v)\, d \mu_0(v) -\int_{M}\langle y, \nabla \phi_0(y)\, d \nu_0(y)
\\
&=&\int_M \langle \nabla {\tilde \phi_1}(y), y\rangle \, d\nu_0(y)-\int_{M}\langle y, \nabla \phi_0(y)\, d \nu_0(y)\\
&=& \int_M \langle \nabla {\tilde \phi_1}(y)-\nabla \phi_0(y), y\rangle \, d\nu_0(y), 
\end{eqnarray*}
where $ {\tilde \phi_1}$ is the concave Legendre transform of $\phi_1$.

\section{Duality for  Ballistic Transports}

We now use the Hopf-Lax formulae established in the previous section to prove Theorem \ref{duality}. The other main ingredient is the duality formula exhibited by Bernard and Buffoni \cite{B-B} for the optimal mass transport $C_T(\nu_0, \nu_T)$, where $\nu_0$ and $\nu_T$ are two given probability measures on $M$. \\

\noindent{\bf Proof of Theorem \ref{duality}:} To prove the duality formula (\ref{ballistic.dual1}), first note that if $V_0$ is any initial Kantorovich potential for $C_T(\nu_0,\nu_T)$, then the final one can be taken to be 
\begin{eqnarray*}
V_T(x)&=&\inf_{y\in M} c_T(y, x)+V_0(y) \\
&=&\inf\Big\{V_0(\gamma (0))+\int_0^TL(s,\gamma (s), {\dot \gamma}(s))\, ds; \gamma \in C^1([0, T), M);   \gamma(T)=x\Big\}.
\end{eqnarray*}
 In other words, $V(T, x)$ is the final state (at time $T$) of a variational solution of (\ref{HJ.0}) starting at $V_0$. 

Now use the Hopf-Lax formula to write 
${\underline B}_T(\mu_0,\nu_T)=C_T(\nu_0, \nu_T)+{\underline W}(\mu_0, \nu_0),$
for some probability measure $\nu_0$ on $M$. The proof of Theorem \ref{interpol} also yields that 
\begin{equation}\label{one.bisbis}
{\underline B}_T(\mu_0,\nu_T)=\int_Mh(x)\, d\nu_T(x)-\int_{M^*}g(v)\, d\mu_0(v),
\end{equation}
where $g$ (resp $h$) is a convex initial (resp., final) Kantorovich potential for ${\underline B}_T(\mu_0,\nu_T)$ if and only if
\[
{\underline C}_T(\nu_0,\nu_T)=\int_Mh(x)\, d\nu_T(x)-\int_{M}{\tilde g}(v)\, d\mu_0(v),
\]
where $\tilde g(y)=- g^*(-y)$ (resp $h$) is a concave initial (resp., final) Kantorovich potential for ${\underline C}_T(\nu_0,\nu_T)$.
Note that $\tilde g$ is the concave Legendre transform of $-g$. In other words, the duality formula (\ref{ballistic.dual1}) can then be obtained by taking as initial Kantorovich potential for ${\underline C}_T(\nu_0,\nu_T)$ any concave functional $V_0$. The corresponding final Kantorovich potential for ${\underline C}_T(\nu_0,\nu_T)$ is then the final state $V_T$ of a variational solution of (\ref{HJ.0}) starting at $V_0$. As for ${\underline B}_T(\mu_0,\nu_T)$, we then have
\begin{equation*}\label{two}
{\underline B}_T(\mu_0,\nu_T)=\sup\left\{\int_MV_T (x)\, d\nu_T(x)+ \int_{M^*}{\tilde V_0}(v)\, d\mu_0(v); \,  \hbox{$V_0$ concave \& $V_t$  solution of (\ref{HJ.0})} \right\}.
\end{equation*}
The proof is similar for the duality formula (\ref{ballistic.dual2}), provided one uses the Hopf-Lax formula (\ref{HL.two}), replace $L$ by $\tilde L$ and note that the corresponding Hamiltonian is now 
$H_{\tilde L}(q,x)=-H_L(x, q)$.

\section{Optimal maps for the ballistic cost}

We have seen in Section 4, that since the cost function $b_T$ is continuous, the infimum ${\underline B}_T(\mu_0,\nu_T)$ is attained at some probability measure $\pi_0\in {\mathcal K}(\mu_0, \nu_T)$, and that the supremum in the dual problem  is attained at some pair $(g, h)$ of admissible Kantorovich functions, where $g$ is convex. In other words, the optimal transport plan $\pi_0$ is supported on the set
\[
{\mathcal O}:=\{(v, x)\in M^*\times M; \, b_T(v, x)=h(x)-g(v)\}
\]
Moreover, for each $(v, x) \in {\mathcal O}$, the function $w\to -g(w)+h(x)-b_T(w,x)$ attains its maximum at $v$, which means that $
-\nabla g(v)\in {\partial}_v b_T(v, x).$ By Proposition \ref{key}, there exists a Hamiltonian trajectory $(\gamma (t), \eta (t))$ over $[0, T]$ starting at $(-\nabla g(v), v)$ and ending at $x$. The rest of the proof of Theorem \ref{attain} is clear. It is however instructive to make the connection with the known results 
regarding $c_T$.

Indeed, if $c_T$ satisfies {\it the twist condition}, that is if there is at most one Hamiltonian trajectory starting at a given initial state $(y, v)$, then $x$ will be determined by $v$, and ${\mathcal O}$ will be supported by a graph. 
This is indeed the case for Tonelli Lagrangians, which were considered
in the compact case by Bernard-Buffoni \cite{B-B}, and by Fathi-Figalli \cite{F-F} in the case of a Finsler manifold. 

\begin{defn} \rm $L$ is said to be a {\it Tonelli Lagrangian on $M\times M$}, if it satisfies the following properties:
 \begin{enumerate}
\item[(a)] $L$ is C$^2$;
\item[(b)] for every $x \in M$, the function $v\to L(x, v)$ is strictly convex on $M$;
\item[(c)] there exist a constant $C>-\infty$ and a non-negative function $\theta: \R^d\to \R$ with superlinear growth, i.e., $\lim\limits_{|v|\to +\infty}\frac{\theta (v)}{|v|}=+\infty$, sich that $L(x,v) \geq C+\theta (v).$
\end{enumerate}
\end{defn}
We also recall the following \cite[Definition 5.5.1, page 129]{A-G-S}:
\begin{defn} 
\label{approxdiff}
{\rm
Say that $f :M \rightarrow \R$ has an \textit{approximate differential} at $x \in M$
if there exists a function $h:M \rightarrow \R$
differentiable at $x$ such that the set $\{f = h\}$ has density $1$ at $x$ with respect to the Lebesgue measure.
In this case, the approximate value of $f$ at $x$ is defined as $\tilde f(x)=h(x)$, and the approximate differential of $f$ at $x$ is defined as $\tilde d_xf=d_xh$.
  It is not difficult to show that this definition makes sense. In fact, both $h(x)$, and $d_xh$ do not depend on the choice of $h$, provided $x$ is a density point of the set $\{f = h\}$.
}
\end{defn}
If $L$ is a Tonelli Lagrangian,
 the Hamiltonian $ H: M\times M^*\to\R$ is then C$^1$, and  the Hamiltonian vector field $X_H$ on $M\times M^*$ is then 
$$
X_H(x,v)=(\frac{\partial H}{\partial v}(x,v),-\frac{\partial H}{\partial x}(x,v)), 
$$
and the associated system of ODEs is given by
\begin{equation}\label{Ham}
\left\{
\begin{array}{l}
\dot x=\dfrac{\partial H}{\partial v}(x,v)\\
\dot v=-\dfrac{\partial H}{\partial x}(x,v).
\end{array}
\right.
\end{equation}
The connection between minimizers $\gamma:[a,b]\to M$ of $I_L$ and solutions of  $(\ref{Ham})$ is as follows. If we write 
$$
x(t)=\gamma(t)\quad \text{and} \quad
v(t)=\frac{\partial L}{\partial p}(\gamma(t),\dot\gamma(t)), 
$$
then 
$x(t)=\gamma(t)$ and $v(t)$ are C$^1$ with $\dot x(t)=\dot\gamma(t)$, and the Euler-Lagrange equation yields 
$
\dot v(t)=\frac{\partial L}{\partial x}(\gamma(t),\dot\gamma(t)),
$
from which follows 
that $t\mapsto (x(t),v(t))$ satisfies (\ref{Ham}). Note also that since $L$ is a Tonelli Lagrangian, the Hamiltonian $H$ is actually C$^2$, and  the vector field $X_H$ is C$^1$. It therefore defines a (partial) C$^1$ flow $\phi^H_t$. 

There is also a (partial) {\rm C}$^1$  flow $\phi^L_t$ on $M\times M^*$ such that every speed curve of an $L$-minimizer is a part of an orbit of $\phi^L_t$. This flow is called the Euler-Lagrange flow, is defined by
$$\phi^L_t=\Leg^{-1}\circ \phi^H_t \circ \Leg,$$
where $
\Leg: M\times M \rightarrow M\times M^*,
$
is the global Legendre transform
$
(x,p) \mapsto (x,\frac{\partial L}{\partial p}(x,p)).
$
Note that $\Leg$ is a homeomorphism on its image whenever $L$ is a Tonelli Lagrangian. We now recall the following result. 

\begin{thm}\label{FF} (Fathi-Figalli \cite{F-F}) {\sl Let $L$ be a Tonelli Lagrangian on $M$.
Fix $T>0$, $\nu_0,\nu_T$ a pair of probability measure on $M$, with $\nu_0$ absolutely continuous with respect to Lebesgue measure.
Then there exists a uniquely $\nu_0$-almost everywhere defined transport map $S:M\to M$ from $\nu_0$ to $\nu_T$
which is optimal for the cost $c_T$. Moreover, any plan $\gamma \in {\mathcal K}(\nu_0,\nu_T)$, which is optimal for the cost $c_T$, verifies $\gamma(\operatorname{Graph}(S))=1$.
\begin{enumerate}
\item If 
$(h, k)$ is an optimal Kantorovich pair, that is if 
\[
h(x)-k(y)=c_T(y,x) \quad \hbox{for $\gamma$-a.e. $(y,x)$ in $M\times M$,}
\]
then there is a Borel set $B$ of full $\nu_0$-measure, such that the approximate differential
$\tilde d_yk$  of $k$ at $y$ is defined for $y\in B$, the map $y\mapsto \tilde d_yk$
is Borel measurable on $B$,
and the transport map $S$ is defined on $B$ (hence $\mu$-almost everywhere) by
$$
S(y)=\pi^*\phi^H_T(y,\tilde d_yk),
$$
where $\pi^*:M\times M^*\to M$ is the canonical projection, and $\phi^H_t$ is the Hamiltonian flow of the Hamiltonian $H$ associated to $L$.

\item A Lagrangian description for $S$ valid on $B$  (hence $\nu_0$-almost everywhere) is given by
$$S(y)=\pi\phi^L_t(y,\widetilde\grad^L_yk),$$
where $\pi=\pi^*\circ\Leg$, $\phi^L_t$ is the Euler-Lagrange flow of $L$, and $y\to\widetilde\grad^L_yk$
is the measurable vector field on $M$ defined on $B$ by
$$
\frac{\partial L}{\partial v}(y,\widetilde\grad^L_yk)=\tilde d_yk.
$$
Moreover, for  every $y\in B$, there is a unique $L$-minimizer $\gamma:[0,T]\to M$, with $\gamma(0)=y, \gamma(T)=S(y)$,
and this curve $\gamma$ is given by $\gamma(s)=\pi\phi^L_s(y,\widetilde\grad^L_yk)$, for $0\leq s\leq T$.
\end{enumerate}
}
\end{thm}
We now give the following corresponding result for $b_T$. 

\begin{thm} In addition to $(A_1)$, assume that $L$ is a Tonelli Lagrangian, then 
\begin{enumerate}
\item There exists a concave function $k: M \to \R$ such that   
\begin{equation}
{\underline B}_T(\mu_0,\nu_T)=\int _{M^*} b_T \large(v, S_T \circ \nabla {\tilde k}(v)\large) d\mu_0(v).
\end{equation}
where 
$$
S_T(y)=\pi^*\phi^H_T(y,\tilde d_y k),
$$
$\pi^*:M\times M^*\to M$ being the canonical projection, ${\tilde k}$ is the concave Legendre transform of $k$, and $\phi^H_t$ the Hamiltonian flow  associated to $L$.

In other words, an optimal map for ${\underline B}_T(\mu_0,\nu_T)$ is given by $v\to \pi^*\phi^H_T(\nabla {\tilde k}(v), v)$.

\item There exists a convex function $h: M^* \to \R$ such that 
\begin{equation}
{\overline B}_T(\mu_0,\nu_T)= \int _{M^*} b_T \large({S}^*_T \circ \nabla h^*(x), x\large) d\nu_T(x),
\end{equation}
where 
${S}^*_T(v)=\pi^*\phi^{H_*}_T(v, \tilde d_vh),$
and $\phi^{H_*}_t$ the flow  associated to the Hamiltonian $H_*(v,x)=-H(-x, v)$, whose Lagrangian is $L_*(v,q)=L^*(-q, v)$. 

In other words, an optimal map for ${\overline B}_T(\mu_0,\nu_T)$ is given by the inverse of the map $x\to \pi^*\phi^{H_*}_T(\nabla h^* (x), x)$.

\item We also have 
\begin{equation}
{\overline B}_T(\mu_0,\nu_T)= \int _{M^*} b_T (v, \nabla h \circ {\tilde S}_Tv\large) d\mu_0(v),
\end{equation}
where 
$$
{\tilde S}_T(v)=\pi^*\phi^{{\tilde H}}_T(v, {\tilde d}_vh_0)
$$ 
$\phi^{{\tilde H}}_t$ being the Hamiltonian flow associated to $\tilde L$ (i.e., ${\tilde H}(v, x)=-H(x, v)$, and $h_0$ the solution $h(0, v)$ of the reverse Hamilton-Jacobi equation (\ref{dHJ}) with $h(T, v)=h(v)$. 

\end{enumerate}
\end{thm}
\noindent{\bf Proof:} Start again by the interpolation inequality and write that 
\[
{\underline B}_T(\mu_0,\nu_T)=C_T(\nu_0, \nu_T)+{\underline W}(\mu_0, \nu_0),
\]
for some probability measure $\nu_0$. The proof also shows that there exists a concave function $k: M \to \R$ and another function $h:M\to \R$ such that $(\nabla {\tilde k})_\#\mu_0=\nu_0$,
\[
{\underline W}(\mu_0, \nu_0)=\int_M \langle \nabla {\tilde k}(v), v\rangle d\mu_0(v).
\]
and
\[
C_T(\nu_0, \nu_T)=\int_Mh(x)\, d\nu_T(x)-\int_M k(y)\, d\nu_0(y).
\]
Now use the theorem of Fathi-Figalli to write 
\begin{equation}
C_T(\nu_0, \nu_T)=\int _{M} c_T (y, S_Ty) d\nu_0(y), \end{equation}
where $S_T(y)=\pi^*\phi^H_T(y,\tilde d_y k)$. 
Note that 
\begin{equation}
\hbox{${\underline B}_T(\mu_0,\nu_T)\leq \int _{M^*} b_T (v, S_T\circ \nabla {\tilde k} (v)) d\mu_0(v),$}
\end{equation}
since $\nabla {\tilde k}_\#\mu_0=\nu_0$ and $(S_T)_\#\nu_0=\nu_T$, and  therefore $(I \times S_T\circ \nabla {\tilde k})_\#\mu_0$ belongs to ${\mathcal K}(\mu_0, \nu_T)$.\\
On the other hand, since $ b_T(v, x)\leq  c_T(\nabla {\tilde k}(v), x )+\langle \nabla {\tilde k}(v), v\rangle$ for every $v\in M^*$, we have
\begin{eqnarray*}
{\underline B}_T(\mu_0,\nu_T)&\leq& \int _{M^*} b_T (v, S_T\circ \nabla {\tilde k} (v)) d\mu_0(v)\\
&\leq& \int_{M^*}\{c_T(\nabla {\tilde k}(v), S_T\circ \nabla {\tilde k} (v)) +\langle \nabla {\tilde k}(v), v\rangle \}\, d \mu_0(v)\\
&=&\int _{M} c_T (y, S_Ty) d\nu_0(y) + \int_{M^*}\langle \nabla {\tilde k} (v), v\rangle \, d \mu_0(v)   \\
&=& C_T(\nu_0, \nu_T)+{\underline W}(\mu_0, \nu_0)\\
&=&{\underline B}_T(\mu_0,\nu_T).
\end{eqnarray*}
It follows that 
\begin{equation*}
{\underline B}_T(\mu_0,\nu_T)=\int _{M^*} b_T (v, S_T\circ \nabla {\tilde k} (v)) d\mu_0(v)=\int _{M^*} b_T (v, \pi^*\phi^H_T(\nabla {\tilde k} (v),\tilde d_{\nabla {\tilde k} (v)}k d\mu_0(v). 
\end{equation*}
Since $k$ is concave, we have that $\tilde d_xk=\nabla k(x)$, hence $\tilde d_{\nabla {\tilde k}(v)}k=\nabla k\circ \nabla {\tilde k}(v)=v$, which yields our claim that 
$$
 {\underline B}_T(\mu_0,\nu_T)=\int _{M^*} b_T \large(v, \pi^*\phi^H_T(\nabla {\tilde k} (v), v)\large) d\mu_0(v).
 $$

2) The proof for ${\overline B}_T(\mu_0,\nu_T)$ is similar and can be done by just reversing the order and proceeding from $\nu_T$ to $\mu_0$. For 3), there exists a convex function $h: M^* \to \R$ such that 
\begin{equation}
{\overline B}_T(\mu_0,\nu_T)= \int _{M^*} b_T (v, \nabla h\circ {\tilde S}_T(v)\large) d\mu_0(v).
\end{equation}
where 
$$
{\tilde S}_T(v)=\pi^*\phi^{{\tilde H}}_T(v, {\tilde d}_vh_0)
$$ 
$\phi^{{\tilde H}}_t$ being the Hamiltonian flow associated to $\tilde L$, and $h_0$ the solution $h(0, v)$ of the reverse Hamilton-Jacobi equation (\ref{dHJ}) with $h(T, v)=h(v)$. 

As above, we start with the interpolation formula
$${\overline B}_T(\mu_0,\nu_T)={\overline W}(\nu_T, \mu_T)-{\tilde C}_T(\mu_0, \mu_T), $$
where $\mu_T$ is a probability measure on $M^*$. The proof yields  
a convex function $h: M^*\to \R$ and another 
function $h_0:M^*\to \R$ as above such that  $(\nabla h^*)_\#\nu_T=\mu_T$
\[
{\overline W}(\nu_T, \mu_T)=\int_{M^*} \langle \nabla^* h(x), x\rangle d\nu_T(x), 
\]
and
\[
{\tilde C}_T(\mu_0, \mu_T)=\int_Mh(w)\, d\mu_T(w)-\int_{M^*}h_0(v)\, d\mu_0(v).
\]
Now use the theorem of Fathi-Figalli to write 
\begin{equation}
{\tilde C}_T(\mu_0, \mu_T)=\int _{M} {\tilde c}_T (v, {\tilde S}_Tv) d\mu_0(v),
\end{equation}
where ${\tilde S}_T(v)=\pi^*\phi^{{\tilde H}}_T(v,\tilde d_vg)$, and 
${\tilde H} (v, x):=H_{\tilde L} (v, x)=-H_L(x,v)$. 

Note  that 
\begin{equation}
\hbox{${\overline B}_T(\mu_0,\nu_T)\geq \int _{M^*} b_T (v, \nabla h\circ {\tilde S}_T(v) ) d\mu_0(v),$}
\end{equation}
since $(\nabla h)_\#\mu_T=\nu_T$ and $({\tilde S}_T)_\#\mu_0=\mu_T$, and  therefore $(I \times \nabla h\circ {\tilde S}_T)_\#\mu_0$ belongs to ${\mathcal K}(\mu_0, \nu_T)$.

On the other hand, since $ b_T(v, y)\geq \langle \nabla h(x), y\rangle -{\tilde c}_T(v, \nabla h(x) )$ for every $(x, y, v)\in M\times M\times M^*$, we have for every $(x, v)\in M\times M^*$, 
$$
 b_T(v,  \nabla h\circ {\tilde S}_T(v))\geq \langle \nabla h(x),  \nabla h\circ {\tilde S}_T(v)\rangle -{\tilde c}_T(v, \nabla h(x) ),
$$
and by taking $x=\nabla h^* \circ {\tilde S}_Tv$, we have 
$$
 b_T(v,  \nabla h\circ {\tilde S}_T(v))\geq \langle {\tilde S}_Tv,  \nabla h\circ {\tilde S}_T(v)\rangle -{\tilde c}_T(v,  {\tilde S}_T(v) ).
$$
It follows that 
\begin{eqnarray*}
{\overline B}_T(\mu_0,\nu_T)&\geq& \int _{M^*} b_T (v, \nabla h\circ {\tilde S}_Tv ) d\mu_0(v)\\
&\geq& \int_{M^*}\{\langle {\tilde S}_Tv,  \nabla h\circ {\tilde S}_Tv\rangle -{\tilde c}_T(v,  {\tilde S}_Tv)\}\, d \mu_0(v)\\
&=&\int_{M^*}\langle v,  \nabla h(v)\rangle \, d\mu_T(v)  -\int_{M^*}{\tilde c}_T(v,  {\tilde S}_Tv)\, d \mu_0(v) \\
&=& {\overline W}(\nu_T, \mu_T)-{\tilde C}_T(\mu_0, \mu_T)\\
&=&{\overline B}_T(\mu_0,\nu_T), 
\end{eqnarray*}
and therefore,
\begin{equation*}
{\underline B}_T(\mu_0,\nu_T)=\int _{M^*} b_T (v, \nabla h\circ {\tilde S}_Tv ) d\mu_0(v)=\int _{M^*} b_T (v, \nabla h\circ \pi^*\phi^{{\tilde H}}_T(v,\tilde d_vh_0) ) d\mu_0(v).
\end{equation*}
\section{The Eulerian formulation}

For the sake of brevity, we shall assume that $L$ is a Tonelli Lagrangian and use the following characterization of $C_T(\nu_0, \nu_T)$ given in \cite{B-B} in the compact setting. See the thesis of Schachter \cite{Sc} for the case $M=\R^d$. 
\begin{eqnarray*}
C_T(\nu_0,\nu_T)
&=& \inf\left\{\int_0^T \int_{M} L\bigl(x, w_t  (x) \bigr) d\varrho_t(x) dt; \, (\varrho, w) \in P(0,T; \nu_0, \nu_T)\right\},  
 \end{eqnarray*}
where $P(0,T; \nu_0, \nu_T)$ is the set of pairs $(\varrho, w)$ such that $t \rightarrow \varrho_t \in \mathcal P(M)$, $t \rightarrow w_t \in \R^n$ are paths of Borel vector fields such that   
\begin{eqnarray}
\partial_t \varrho+ \nabla \cdot (\varrho w)&=&0 \quad {\rm in}\, \mathcal D'\bigl((0,T) \times M \bigr)\\\varrho_0=\nu_0&,&\varrho_T=\nu_T.\nonumber
\end{eqnarray} 
Combining this with (\ref{HL.two}), we get that 
\begin{eqnarray*}
{\underline B}_T(\mu_0,\nu_T)&=&
\inf\{{\underline W}(\mu_0, \nu})+ C_T({\nu, \nu_T);\, \nu\in {\mathcal P}(M)\}\\
&=& \inf\left\{ {\underline W}(\mu_0, \nu) +
\int_0^T \int_{M} L\bigl(x, w_t  (x) \bigr) d\varrho_t(x) dt; \,  \nu\in {\mathcal P}(M),\, (\varrho, w) \in P(0,T; \nu, \nu_T)\right\}\\  
&=& \inf\left\{ {\underline W}(\mu_0, \rho_0) +
\int_0^T \int_{M} L\bigl(x, w_t  (x) \bigr) d\varrho_t(x) dt; \, (\varrho, w) \in P_{ter}(0,T; \nu_T)\right\},  
 \end{eqnarray*}
where $P_{ter}(0,T; \nu_T)$ is the set of pairs $(\varrho, w)$ such that 
\begin{eqnarray}
\partial_t \varrho+ \nabla \cdot (\varrho w)&=&0 \quad {\rm in}\, \mathcal D'\bigl((0,T) \times M \bigr)\\\varrho_T&=&\nu_T.\nonumber
\end{eqnarray} 
The same reasoning holds for ${\overline B}_T(\mu_0,\nu_T)$.

\end{document}